\RequirePackage{fix-cm}
\documentclass[smallextended]{svjour3}       
\smartqed  
\usepackage{graphicx}

\usepackage{amsmath}
\usepackage{amssymb}
\usepackage{latexsym}
\usepackage{psfrag}
\usepackage{subfigure}
\usepackage{bm}

\def\eqref#1{(\ref{#1})}

\def\bd{\mathbf{d}}

\def\bbf{\mathbf{f}}
\def\bg{\mathbf{g}}
\def\bh{\mathbf{h}}

\def\bn{\mathbf{n}}

\def\br{\mathbf{r}}

\def\bu{\mathbf{u}}
\def\bv{\mathbf{v}}
\def\bw{\mathbf{w}}
\def\bx{\mathbf{x}}

\def\bA{\mathbf{A}}
\def\bB{\mathbf{B}}

\def\bF{\mathbf{F}}

\def\bP{\mathbf{P}}

\def\bR{\mathbf{R}}

\def\bxi{\bm \xi}

\def\n{\noindent}
\def\pt{\partial}

\def\Re{\mathbb{R}}

\def\minmod{\mbox{minmod}}
\def\f#1#2{\frac {#1}{#2}}
\def\f32{\frac 32}

\def\d{\displaystyle}

\def\b#1{\overline#1}
\def\bga{\begin{array}}
\def\eda{\end{array}}

\def\De{\Delta}

\def\ev{\equiv}
\def\ol{\overline}

\def\gm{\gamma}

\def\th{\theta}
\def\nb{\nabla}
\def\ep{\epsilon}
\def\la{\lambda}
\def\La{\Lambda}
\def\sg{\sigma}

\def\al{\alpha}
\def\td{\tilde}
\def\rw{\rightarrow}
\def\iy{\infty}
\def\Om{\Omega}
\def\om{\omega}

\def\sq{\sqrt}

\def\d{\displaystyle}
\def\dfr#1#2{\displaystyle{\frac{#1}{#2}}}

 \newtheorem{thm}{Theorem}[section]

 \newtheorem{prop}[thm]{Proposition}
 \newtheorem{rem}[thm]{Remark}

\begin{document}

\title{Two-Stage Fourth Order: Temporal-Spatial Coupling in Computational Fluid Dynamics (CFD)
}


\author{Jiequan Li
}


\institute{Jiequan Li\at
            Laboratory of Computational Physics,  Institute of Applied Physics and Computational Mathematics, Beijing, P. R. China; Center for Applied Physics and Technology, Peking University, Beijing, P. R. China  \\
              Tel.: +86-10-61935465\\
              \email{li\_jiequan@iapcm.ac.cn}           
}

\date{Received: date / Accepted: date}

\maketitle
\markboth{Two-stage fourth order: Temporal Spatial Coupling in CFD}{Jiequan Li}

\begin{abstract}
With increasing engineering demands, there need high order accurate schemes embedded with precise physical information in order to capture delicate small scale structures and strong waves with correct``physics". There are two families of high order methods: One is the method of line, relying on the Runge-Kutta (R-K) time-stepping. The building block is the Riemann solution labeled as the solution element ``$1$". Each step in R-K just has first order accuracy. In order to derive a fourth order accuracy scheme in time, one needs four stages labeled as ``$1\odot 1\odot 1\odot 1=4$". The other is the one-stage Lax-Wendroff (LW) type method, which is more compact but is complicated to design numerical fluxes and hard to use when applied to highly nonlinear problems.
In recent years, the pair of solution element and dynamics element, labeled as ``2", are taken as the building block. The direct adoption of the dynamics implies the inherent temporal-spatial coupling. With this type of building blocks, a family of two-stage fourth order accurate schemes, labeled as ``$2\odot2=4$", are designed for the computation of compressible fluid flows. The resulting schemes are compact, robust and efficient. This paper contributes to elucidate how and why high order accurate schemes should be so designed. To some extent, the ``$2\odot 2=4$" algorithm extracts the advantages of the method of line and one-stage LW method. As a core part, the pair ``$2$" is expounded and LW solver is revisited. The generalized Riemann problem (GRP) solver, as the discontinuous and nonlinear version of LW flow solver, and the gas kinetic scheme (GKS) solver, the microscopic LW solver, are all reviewed. The compact Hermite-type data reconstruction and high order approximation of boundary conditions are proposed. Besides, the computational performance and prospective discussions are presented.
\end{abstract}

\section{Introduction}
\label{sec:intro}
In  the simulation of compressible fluid flows or related problems, there are two families of commonly-used high order accurate numerical schemes:  One is the family of  methods of line, for which the fluid dynamical system is written in semi-discrete form and  the Runge-Kutta (RK) temporal iteration is employed for the temporal discretization, such as RK-WENO \cite{Jiang}, RK-DG \cite{DG-2} and their variants.  The building blocks comprises of the solution element,  the associated Riemann solution, which is labeled as ``$1$" in order to be in contrast with the Lax-Wendroff (LW) type flow solvers. The fourth order RK temporal iteration is labeled as ``$1\odot 1\odot 1\odot1 =4$. This family of schemes have very favorable properties such as simplicity in time-stepping  for complex engineering problems.  The limitation is also obvious such as  compactness, efficiency and fidelity.  The other is the family of one-stage LW type methods, the numerical realization of Cauchy-Kowalevski (CK) approach \cite{Evans-2002} for the corresponding partial differential equations.  This family of methods have the strong temporal-spatial coupling property, leading to very compact numerical schemes. However, when applied to high nonlinear problems, the complex construction of numerical fluxes hampers to develop high (more than two) order accurate schemes. Particularly, as strong waves (discontinuities ) are present in flows (solutions), the  CK procedure loses its physical and mathematical meanings,  exhibiting  the instability of the resulting schemes near discontinuities.
\vspace{0.2cm}

Careful inspection of these two families of methods motivates to combine the merits of both methods:  The simplicity of multi-stage RK methods and the temporal-spatial coupling of LW type methods.  This straightforward combination immediately yields a two-stage fourth order accurate temporal discretization for the LW type flow solvers \cite{Li-Du-2016}, which is labeled as "$2\odot 2=4$".  Here "$4$"  just represents "fourth" order accurate temporal discretization, but ``$2$" has much deeper implications, some of which are enumerated below.

\begin{enumerate}
\item[(i)] {\bf ``2" represents a pair.}  Unlike the methods of line, this method adopts the pair, the conservative variables and their dynamics, e.g., the velocity and the acceleration, as the building block to design numerical  schemes. In \cite{Li-Du-2016}, we call this pair as the Riemann solver and the LW type solver.
\vspace{0.2cm}

\item[(ii)] {\bf ``2" implies the temporal-spatial coupling.}  The LW flow solver implies the temporal-spatial coupling property of resulting schemes. This is necessary to simulate the temporal-spatial coherent structures of fluid flows.
\vspace{0.2cm}

\item[(iii)] {\bf ``2" stands for second order accuracy in time.} Of course,``2" also  symbolizes the temporal accuracy of resulting schemes and requires that at each of the two stages the above pair should be the building block.
\vspace{0.2cm}

\item[(iv)] {\bf ``2" indicates the exchange of kinematics and thermodynamics.   }  The Gibbs relation plays a fundamental role in compressible fluid flows. In the dynamical process, there is always the interaction of kinematics and thermodynamics.  The stronger nonlinear waves, e.g. shocks, exist in the fluid flows, the more fundamental  role the thermodynamics plays.

\item[(v)] {\bf ``2" guarantees the compactness and efficiency. }  Since only two stages are taken to achieve fourth order temporal accuracy, half amount of spatial discretization treatments are saved and much smaller computational stencils are needed. Hence the  resulting schemes are more compact and   efficient.

\item[(vi)] {\bf ``2" reflects the consistency of mathematical and physical expressions of fluid dynamics.}  The fluid dynamical systems  essentially consist of balance laws,  which say the relation between the change rate of physical quantities and the associated fluxes. The form of balance laws  always makes sense no matter whether there are discontinuities in the solution.  The Lax-Wendroff type flow solvers inherently reflects the consistency between the  physical implication of fluid dynamical systems and their mathematical formulation.
\end{enumerate}

In this paper we will elucidate the idea of this new family of schemes by interpreting the philosophy from ordinary differential equations (ODEs) to fluid dynamical systems,  reviewing the well-used GRP and GKS solvers as  the representatives of the Lax-Wendroff  type solvers, building high order temporal-spatially coupled high order accurate schemes with favorable computational performance.

We organize this paper in the following sections.  In Section 2, we propose this new family of methods and the corresponding ``$2\odot 2$" algorithm.  In Section 3,  we review the generalized Riemann problem (GRP) solver and in Section 4 continue to review a kinetic solver, the gas kinetic scheme (GKS) solver. In Section 5,  we  introduce the compact Hermite-type interpolation for the data reconstruction. In Section 6, we discuss the approximation of boundary conditions to suit for the $2\odot 2$ algorithm. In Section 7, we remark the computational performance of this approach in terms of computational efficiency, robustness and fidelity.

\section{What is "$2\odot 2=4$"? }
This section serves to elucidate the meaning of "$2\odot 2=4$" for hyperbolic problems and particularly compressible fluid flows, and review the two-stage fourth order accurate schemes proposed in \cite{Li-Du-2016}.  We remind that this strategy may not be suitable for incompressible flows or it needs some modifications but certainly awaits for further improvement.

\subsection{Start with ODEs and philosophic thinking}
Let's recall the Runge-Kutta (RK) method for an  ordinary differential equation
\begin{equation}
\dfr{dy}{dt} =f(t,y).
\label{eq:ode}
\end{equation}
The Runge-Kutta method takes the iteration procedure
\begin{equation}
\begin{array}{l}
\d y_{n+1} =y_n + h\sum_{i=1}^s b_ik_i\\[3mm]
\d k_i = f\left(t_n+c_i h, y_n +h\sum_{j=1}^{i-1} a_{ij}k_j\right), \ \ i=1,\cdots, s,
\end{array}
\label{eq:RK}
\end{equation}
where $h$ is the time increment, $a_{ij}$, $b_i$ and $c_i$ satisfy the Butcher tableau \cite{Butcher-1964}. The building block of RK is the solution element $y$.  In order to devise a $s$-th order accurate scheme, one needs s-stage iteration, which is parameter-dependent.  In this paper, we focus on fourth order accurate schemes and therefore label the fourth order RK schemes as $1\odot 1\odot 1\odot 1=4$.  The notation ``$\odot$" is an operation satisfying certain requirement such as stability.
\vspace{0.2cm}

The RK method lays the foundation of numerical approximations to ODEs. Note that this method only uses the solution element ``$y$", but ignores the dynamics element $dy/dt$. This sounds confusing,  however, one may pay his attention to the role of the dynamics element if he is familiar with  the symplectic algorithm for Hamitonian system \cite{Feng-2010} for which   the pair of the position and momentum are together used for the computation in order to preserve the symplectic structure. The momentum can be regarded as the dynamics element of the position (solution element). The word ``{\em symplectic} " itself has the meaning of ``{\em pair}".
\vspace{0.2cm}

With this philosophical thinking, it is reasonable to construct multi-stage multi-derivative algorithms for dynamical systems (ODEs). Indeed, this was achieved in \cite{Tsai-2010} with many subsequent works \cite{Tsai-2014,Seal-2014,C-Seal-2016,Seal-2018}. The  building block is the pair, the solution element  and the dynamics element (the derivatives).  Specifically, a multi-stage two-derivative algorithm is written as
\begin{equation}
\begin{array}{l}
\d Y_i =y_n +h\sum_{j=1}^{i-1} a_{ij} f(Y_j) +h^2 \sum_{j=1}^{i-1} \hat a_{ij} g(Y_j),  \ \ i=1,\cdots, s,\\[3mm]
\d y_{n+1}=y_n + h\sum_{i=1}^s b_i f(Y_i) +h^2 \sum_{i=1}^s \hat{b}_i g(Y_i),
\end{array}
\label{eq:MSMT}
\end{equation}
where we suppress the dependence of $f$ on $t$ for simplicity so that $f=f(y)$, the coefficients $a_{ij}$, $\hat{a}_{ij}$, $b_i$, and $\hat{b}_i$ can be displayed in an extended Butcher tableau \cite{Tsai-2010}. Here  the function  $g(y)=f'(y) f(y)$ is given using  the chain rule
\begin{equation}
g(y) =\dfr{d}{dt}f(y) = f'(y)\dfr{dy}{dt}=f'(y)f(y).
\end{equation}
The dynamical element is implicitly used in the construction of algorithm \eqref{eq:MSMT}.
This is why this method is of multi-stage two-derivative type with the pair $(y,dy/dt)$ as the building block. In particular, as $s=2$, we have the two-stage fourth order accurate time-stepping algorithm in the  form
\begin{equation}
\begin{array}{l}
y_* =y_n +\dfr{h}2 f(y_n) +\dfr{h^2}8 f'(y_n) f(y_n), \\[3mm]
y_{n+1} =y_n+ h f(y_n) +\dfr{h^2}6\left(f'(y_n)f(y_n) +2f'(y_*) f(y_*)\right),
\end{array}
\label{eq:2+2-scheme}
\end{equation}
 labelled as the ``$2\odot 2=4$" algorithm, which was independently derived in \cite{Li-Du-2016} for hyperbolic conservation laws. See the discussion in the subsequent sections. For \eqref{eq:2+2-scheme}, the first ``$2$" represents the two-stage approach, the second ``$2$" means the pair of the solution element and the dynamics element, and ``$4$" stands for the fourth order accurate approximation to \eqref{eq:ode}.  Certainly, the first ``$2$" has more implications when applied to the fluid dynamical systems for compressible flows.  Besides, the notation ``$\odot$" is used here to symbolize the mathematical operation currently. Probably in the future, this notation could be replaced by a better one.

 \subsection{Lax-Wendroff flow solvers}

 The Lax-Wendroff method \cite{LW-1960} plays a fundamental role in the development of high order accurate schemes for hyperbolic equations.  The corresponding scheme is unique in the class of three-point  schemes of second order both in space and time.  The feature of uniqueness implies that it is the reference of high order accurate schemes, and the three-point stencil hints at the compactness. Here we are going to show more fundamentals of this  method, which is taken as the building block or higher order accurate schemes.   Moreover, we would like to present as many details as possible because it is unusual that it has not been received  ``enough" attention  since its birth. Part of the reason may be the presence of oscillations near discontinuities when used to simulated compressible fluid flows although it was modified,   e.g., the flux limiter methods in 1980s (\cite{Harten-1983} and its successors), to be suited for the  capture of discontinuities; part of the reason is, more possibly, the seeming   complexity compared to methods of line.  Even more seriously, the misuse in various contexts, such as diffusion equations and (dispersive) KdV type equations,  leads to many controversial issues.

 \subsection{The revisit of Lax-Wendroff method}  Let us first recall the Lax-Wendroff method \cite{LW-1960}.  Consider the advection equation
 \begin{equation}
 u_t +a u_x=0, \ \  \  t>0,\ \ \  x\in [0,1],
 \label{eq:adv}
 \end{equation}
 where $a$ is a constant. The boundary condition remains to be discussed in Section \ref{sec-BC}.  We   approximate \eqref{eq:adv} by assuming that the solution is sufficiently regular, and take the Taylor series expansion   at any point $(x,t)$ to obtain,
 \begin{equation}
 u(x,t+\De t) =u(x,t)+\De t \dfr{\pt u}{\pt t}(x,t) +\dfr{\De t^2}2\dfr{\pt^2 u}{\pt t^2}(x,t)+\mathcal{O}(\De t^3).
 \label{eq:Taylor}
 \end{equation}
 A key step  is the {\em temporal-spatial  coupling}  technique by  taking use of \eqref{eq:adv} to quantify the differentiation relation between the change rate of $u$ and the spatial variation,
 \begin{equation}
 \dfr{\pt u}{\pt t} = -a \dfr{\pt u}{\pt x},  \ \ \ \   \dfr{\pt^2 u}{\pt t^2} = a^2 \dfr{\pt^2 u}{\pt x^2}.
 \label{eq:wave}
 \end{equation}
 Ignoring  truncation errors of order more than three, the Lax-Wendroff scheme is derived  as (cf. \cite{LW-1960}),
 \begin{equation}
 u_j^{n+1} =u_j^n -\dfr{\la}2(u_{j+1}^n-u_{j-1}^n) +\dfr{\la^2}2 (u_{j+1}^n-2u_j^n+u_{j-1}^n),\ \ \la=a\dfr{\De t}{\De x},
 \label{eq:LW}
 \end{equation}
 where central difference approximations are made to guarantee the spatial accuracy, $u_j^n$ represents the point value $u(x_j,t_n)$  at the grid point $(x_j,t_n)$, $x_j=j\De x$, $t_n=n\De t$, with the spatial and temporal increments $\De x$ and $\De t$.  The Taylor expansion process is the same as that in the Cauchy-Kowaleveski approach (see \cite{Evans-2002}), and therefore \eqref{eq:LW} is regarded as the numerical realization of the CK approach.  Note that this process determines the feature of this scheme, implying its  application only in the range of  hyperbolic problems (local behavior or finite propagation property). Any other extension needs serious and cautious treatments.
 \vspace{0.2cm}

The Taylor expansion also relies on the smoothness of the solution.  The successive differentiation \eqref{eq:adv}  gives rise to the risk in the following sense.

 \begin{enumerate}
 \item[(i)] Once  equation \eqref{eq:adv} admits discontinuities in the solution,  the manipulation for \eqref{eq:wave} does not make any sense.  This is the main reason that \eqref{eq:LW} produces oscillations near discontinuities \cite{LW-1960}.
 \vspace{0.2cm}

 \item[(ii)] As this method is applied to highly nonlinear dynamical systems, this manipulation becomes horrible and hampers to develop higher order accurate schemes, due to the successive differentiations.
 \end{enumerate}

 We will  comment on this manipulation appropriately at later sections. Rather now, we reinspect \eqref{eq:adv} and  \eqref{eq:Taylor} from another point of view (after ignoring high order truncation errors),  actually in the finite volume framework,
 \begin{equation}
 \begin{array}{rl}
 u(x,t+\De t)  & =u(x,t) +\De t \dfr{\pt}{\pt t}\left[u+\dfr{\De t}2\dfr{\pt u}{\pt t}\right]\\[3mm]
 & =u(x,t) -a\De t \dfr{\pt}{\pt x}\left[u+\dfr{\De t}2\dfr{\pt u}{\pt t}\right],
 \end{array}
 \label{eq:insp}
 \end{equation}
 where the differentiation relation $\frac{\pt }{\pt t} =-a\frac{\pt}{\pt x}$ is applied.  We immediately realize that for any $(x,t)$
 \begin{equation}
 u(x,t) +\dfr{\De t}{2}\dfr{\pt u}{\pt t} =u\left(x,t+\dfr{\De t}2\right)+\mathcal{O}(\De t^2),
 \label{eq:flux-flow}
 \end{equation}
 as long as the solution is smooth in $t$ (temporal direction or flow direction).  Viewing \eqref{eq:insp} in the finite volume framework,  we obtain over the control volume $[x_{j-\frac 12}, x_{j+\frac 12}]\times [t_n,t_{n+1})$, $x_{j+\frac 12}=\frac 12(x_j+ x_{j+1})$,
 \begin{equation}
 \begin{array}{l}
 u_{j+\frac 12}^* :=u(x_{j+\frac 12},t_n)+\dfr{\De t}2 \dfr{\pt u}{\pt t}(x_{j+\frac 12},t_n),\\[3mm]
 u_j^{n+1}=u_j^n -\la (u_{j+\frac 12}^*-u_{j-\frac 12}^*).
 \end{array}
 \label{eq:LW-2}
 \end{equation}
 The prediction of the value $u_{j+\frac 12}^*$ depends on the approximations to $u(x_{j+\frac 12},t_n)$ and  $\frac{\pt u}{\pt t}(x_{j+\frac 12},t_n)$. This is achieved by the {\em Lax-Wendroff solver}.
 \vspace{0.2cm}

 \n{\bf Lax-Wendroff solver.}  {\em  A Lax-Wendroff solver for \eqref{eq:adv} is the numerical algorithm approximating the values
 \begin{equation}
 u_{j+\frac 12}^n: =\lim_{t\rw t_n+0} u(x_{j+\frac 12}, t),  \ \ \ \ \left(\dfr{\pt u}{\pt t}\right)_{j+\frac 12}^n: =\lim_{t\rw t_n+0} \dfr{\pt u}{\pt t} (x_{j+\frac 12}, t)
 \label{eq:pair}
 \end{equation}
 for the given initial data at $t=t_n$ for \eqref{eq:adv}.
 }
 \vspace{0.2cm}

This pair of values actually provide all quite detailed information along the interface $x=x_{j+\frac 12}$  of control volume  and also the flux
\begin{equation}
\begin{array}{rl}
\dfr{1}{\De t} \int_{t_n}^{t_{n+1}} au(x_{j+\frac 12},t)dt &\d  =au\left(x_{j+\frac 12},t_n+\frac{\De t}2\right)+\mathcal{O}(\De t^2)\\[3mm]
& = a\left[u_{j+\frac 12}^n  +\dfr{\De t}2  \left(\dfr{\pt u}{\pt t}\right)_{j+\frac 12}^n \right]+\mathcal{O}(\De t^2).
\end{array}
\end{equation}
 The two formulae \eqref{eq:LW} and \eqref{eq:LW-2} are equivalent for smooth flows. However, the new formulation \eqref{eq:LW-2} is fundamentally different from \eqref{eq:LW} in the following sense.
 \begin{enumerate}
 \item[(i)]  The formulae \eqref{eq:LW-2} is actually the finite volume formulation for \eqref{eq:adv}. The formulation is more straightforward for fluid dynamical systems than other formulations because it is just the numerical version of balance laws and allows discontinuities as its solution.

 \item[(ii)]  The manipulation \eqref{eq:flux-flow} is legal because the flow should be smooth in time (but not in space), unlike the difference approximation for  LW approach.

 \item[(iii)] The temporal-spatial coupling feature again plays an important role,  e.g.,
 \begin{equation}
 \dfr{\pt u}{\pt t}(x_{j+\frac 12},t_n) =-a\dfr{\pt u}{\pt x}(x_{j+\frac 12}, t_n).
 \label{eq:LW-a}
 \end{equation}
 This feature is crucial for a numerical scheme  to preserve the fluid dynamical properties such as the Galilean invariance.

 \item[(iv)]
 The successive  differentiation \eqref{eq:wave} can be avoided, which is extremely important for nonlinear problems when discontinuities are involved because the manipulation \eqref{eq:wave} makes no sense   both mathematically and physically.
 \end{enumerate}

 For \eqref{eq:adv}, we label ``$2$" for the pair $(u, \frac{\pt u}{\pt t})$ in the Lax-Wendroff solver \eqref{eq:pair}, which is  the building block, as we see,  in the Lax-Wendorff scheme.  It is interesting to observe that \eqref{eq:pair} can be approximated in an  upwind or central way.  The upwind approximation can avoid superfluous information in the scheme.

 \subsection{Lax-Wendroff flow solvers for nonlinear hyperbolic balance laws} We  consider hyperbolic conservation laws
 \begin{equation}
 \bu_t+\bbf(\bu)_x=0,
 \label{eq:HYP}
 \end{equation}
 where the vector $\bu$ is the conservative variable. The natural formulation of \eqref{eq:HYP} is in the finite volume framework,  the balance law over any interval $I_j=(x_{j-\frac 12},x_{j+\frac 12})$,
 \begin{equation}
  \begin{array}{l}
 \dfr{d}{dt}\bar \bu_j(t) =-\dfr{1}{\De x} \left( \bbf(\bu(x_{j+\frac 12},t)) -\bbf(\bu(x_{j-\frac 12},t))\right), \\[3mm]
 \bar\bu_j(t)= \dfr{1}{\De x}\int_{I_j} \bu(x,t)dx,
 \end{array}
 \end{equation}
 or the control volume $(x_{j-\frac 12},x_{j+\frac 12})\times (t_n,t_{n+1})$,
 \begin{equation}
\bar \bu_j^{n+1} =\bar\bu_j^n -\dfr{\De t}{\De x}(\bbf_{j+\frac 12}(t_n; t_{n+1})-\bbf_{j-\frac 12}(t_n; t_{n+1})),
 \end{equation}
 with
 \begin{equation}
  \bar\bu_j^n= \dfr{1}{\De x}\int_{I_j} \bu(x,t_n)dx, \ \ \ \bbf_{j+\frac 12}(t_n; t_{n+1})=\dfr{1}{\De t}\int_{t_n}^{t_{n+1}} \bbf(\bu(x_{j+\frac 12},t))dt.
 \end{equation}
 If one would prefer to other formulations, such as the discontinuous Galerkin method \cite{DG-2}, the following statements still hold.
 \vspace{0.2cm}

 We shift $(x_j,t_n)$ to $(0,0)$ due to the invariance of \eqref{eq:HYP} with respect to the translation of coordinates. In order to proceed in  one of those frameworks,  we have to solve \eqref{eq:HYP}  approximately subject to the initial data
 \begin{equation}
 \bu(x,0) = \bP_\pm(x),  \ \ \mbox{for } \pm x>0,
 \label{data:1D}
 \end{equation}
 where $\bP_\pm(x)$ are smooth functions, typically polynomials, with a jump at $x=0$.  The same as in the linear case \eqref{eq:pair}, a {\em Lax-Wendroff flow solver} for such a problem is an algorithm approximating
 \begin{equation}
 \bu_0:=\lim_{t\rw 0^+}\bu(0,t), \ \ \  \left(\dfr{\pt \bu}{\pt t}\right)_0 =\lim_{t\rw 0^+} \dfr{\pt \bu}{\pt t}(0,t).
 \label{eq:pair-3}
 \end{equation}

 \vspace{0.2cm}

 In general, we consider hyperbolic balance laws in multi-dimensions,
 \begin{equation}
 \bu_t + \nb\cdot \bF(\bu)=\bh, \ \ \ \bF=(\bbf_x, \bbf_y, \bbf_z).
 \label{eq:MHYP}
 \end{equation}
 where $\bh$ is the source term resulting from physics or geometry, $\bx=(x,y,z)$ is the spatial coordinate.   The initial data for \eqref{eq:MHYP} is set to be
  \begin{equation}
 \bu(x,y,z,0) = \bP_\pm(x,y,z), \ \ \ \mbox{for }   \pm \bm\mu \cdot \bx>0,
 \end{equation}
 where $\bm\mu$ is the unit normal of a line or plane $L: \bm\mu \cdot \bx =0$  pointing from the negative side to the positive side, corresponding to the outer normal of  interfaces of computational volume.  The Lax-Wendroff solver for \eqref{eq:MHYP} is to find the pair of values with the same form as in \eqref{eq:pair-3},
 \begin{equation}
 \bu_{L, 0} := \lim_{t\rw 0^+}\bu(L, t), \ \ \  \left(\dfr{\pt \bu}{\pt t}\right)_{L, 0} := \lim_{t\rw 0^+}\dfr{\pt \bu}{\pt t}(L, t),
 \label{eq:pair-2d}
 \end{equation}
 where the limit is taken along the spatial-temporal interface $L\times (0,\De t)$.
 \vspace{0.2cm}

 We want to remark here that the pair $(\bu_{L,0}, (\pt \bu/\pt t)_{L,0})$ can be modulated to any direction in order to suit for an arbitrary Lagrangian-Eulerian  (ALE) method.

 \subsection{Rough comments on the correlation between LW solver and temporal-spatial coupling} The instantaneous temporal derivatives in \eqref{eq:pair-3} and \eqref{eq:pair-2d} can be roughly using the Lax-Wendroff approach
 \begin{equation}
 \dfr{\pt \bu}{\pt t} =-\nb \cdot \bF(\bu) +\bh,
 \end{equation}
 and then $-\nb\cdot \bF(\bu)$ and $\bh$ are approximated using certain technologies such as WENO etc. The same as in \eqref{eq:LW-a}, the coherent relation of  spatial and temporal variations is rooted in this formula.
 \vspace{0.2cm}

 The intuitive outcome of this coupling is  the following.
 \begin{enumerate}
 \item[(i)]  {\em The multidimensional  effect, in particular the transversal effect, is input into  the flux directly.} Thinking of a single advection problem
\begin{equation}
u_t + au_x +bu_y=0,
\end{equation}
where $a, b$ are constants.  For an interface with the normal in the $x$-direction, the transversal effect, expressed in the $y$-direction, is ignored for the standard Riemann solver.  This is further verified for the  wave system
\begin{equation}
p_t+c_0u_x+c_0v_y=0, \ \ \ u_t +c_0p_x =0, \ \ \ v_t +c_0p_y=0,
\end{equation}
where $c_0$ is a constant.  In Table 1, we use three methods to simulate the periodic wave problem. It is observed  that even with the same convergence rate, the RK method produces also ten times of errors than what the second order GRP does for which the transversal effect is included. The solution cannot even converge with the refinement of meshes if only normal flux is used but the transversal effect is not included.  See \cite{Lei-Li-AAM-2018}.

	\begin{table}[!htb]
		\label{tab:wave}
		\centering
		\caption{ $L_1$ error and convergence order of $u$ for the periodic wave problem at  final time $t=2$. with the methods \textit{GRP2D},  \textit{RK2} and \textit{GRP1D}}
	\noindent\begin{tabular}{ccr|cr|cr}
		\hline
		\raisebox{-2.00ex}[0cm][0cm]{$N$} &\multicolumn{2}{c}{\textit{GRP2D}} & \multicolumn{2}{c}{\textit{RK2}} & \multicolumn{2}{c}{\textit{GRP1D}}\\ \cline{2-7}
		 & \multicolumn{1}{c}{$L_1$ error}& \multicolumn{1}{c|}{order}&\multicolumn{1}{c}{$L_1$ error}& \multicolumn{1}{c|}{order}& \multicolumn{1}{c}{$L_1$ error}&\multicolumn{1}{c}{order} \\ \hline	
		$40$ & 4.54E-2	&		 &	1.38E-1	&		 &	4.16E-1	&	   \\[-1mm]
		$80$ & 7.32E-3	&$ 2.63 $&	3.55E-2	&$ 1.96 $&	2.25E-1	&$0.89$\\[-1mm]
		$160$ & 1.33E-3	&$ 2.46 $&	8.96E-3	&$ 1.99 $&	1.17E-1	&$0.94$\\[-1mm]
		$320$ & 2.81E-4	&$ 2.25 $&	2.25E-3	&$ 2.00 $&	6.66E-2	&$0.81$\\[-1mm]
		$640$ & 6.53E-5	&$ 2.10 $&	5.62E-4	&$ 2.00 $&	1.88E-1	&$-1.50$\\[-1mm]
		\hline
			\end{tabular}
			 \end{table}

 \item[(ii)]  {\em The source effect $\bh$ is also reflected through such a process.}   It is simple to see that
 \begin{equation}
 \dfr{\pt u}{\pt t} =-\dfr{\pt f(u)}{\pt x}+h(u,x),
 \end{equation}
 for hyperbolic balance law
  \begin{equation}
 \dfr{\pt u}{\pt t} +\dfr{\pt f(u)}{\pt x}=h(u,x).
 \end{equation}
 This input is essential and indispensable for the well-balancedness, as verified for the shallow water equations \cite{Li-Chen-2006}.
 \end{enumerate}

 There are more fundamentals, such as the thermodynamical effect \cite{Li-Wang-2017}, resulting from the temporal-spatial coupling.

 \subsection{$2\odot 2=4$: Two-stage fourth order accurate schemes}

 In \cite{Li-Du-2016}, the fourth order accurate method is developed for hyperbolic conservation laws. We start with the review of the dynamical system
 \begin{equation}
 \dfr{d}{dt}\bw =\mathcal{L}(\bw),
 \label{eq:dy}
 \end{equation}
 where $\mathcal{L}$ is a linear or nonlinear  operator of $\bw$.  Then we have the following two-stage algorithm for \eqref{eq:dy}.
 \begin{enumerate}
 \item[\bf Stage 1.] Define intermediate values
\begin{equation}
\begin{array}{l}
\bw^* =  \bw^n +\frac 12\De t \mathcal{L}(\bw^n) +\dfr 18 \De t^2 \dfr{\pt}{\pt t}\mathcal{L}(\bw^n),\\
\dfr{\pt}{\pt t}\mathcal{L}(\bw^*) =\dfr{\pt}{\pt \bu}\mathcal{L}(\bw^*) \mathcal{L}(\bw^*),
\end{array}
\label{a-1}
\end{equation}
where the second equation follows from \eqref{eq:dy}, using the chain rule.
\vspace{0.2cm}

\item[\bf Stage 2.]  Advance the solution using the formula
\begin{equation}
\bw^{n+1} =\bw^n + \De t \mathcal{L}(\bw^n) + \dfr 16 \De t^2 \left(\dfr{\pt }{\pt t}\mathcal{L}(\bw^n)+2\dfr{\pt }{\pt t}\mathcal{L}(\bw^*)\right).
\label{a-2}
\end{equation}
\end{enumerate}
This algorithm provides a fourth order accurate approximation to $\bw$.  Originally, this algorithm was proposed in \cite{Tsai-2010,Seal-2016}, and independently in \cite{Li-Du-2016} based  on Lax-Wendroff flow solvers.  Along this direction, one can derive as high order accurate approximations as what he likes \cite{Tsai-2010,Xu-Li-2017}.

\vspace{0.2cm}

When applied to hyperbolic problems \eqref{eq:HYP} and \eqref{eq:MHYP}, one can formulate them in any appropriate framework such as finite volume framework  \cite{Li-Du-2016} or discontinuous Galerkin (DG) framework \cite{Cheng-Li-2018}.  Hence we assume that the computational domain $\Om$ is meshed as $\Om=\d \cup_{j\in\mathcal{J}} \Om_j$ and formulate the problem in the form
\begin{equation}
\dfr{d}{dt} \bw_j(t) =\mathcal{L}_j(\bw),  \ \ \ \ \bw_j =\dfr{1}{|\Om_j|}\int_{\Om_j} \bu(\bx,t)d\bx, \ \ \  \bw=\{\bw_j; j\in \mathcal{J}\}.
\label{eq:semi}
\end{equation}
Thus, this problem boils down to the dynamical system in the form \eqref{eq:dy}. Then we have a two-stage fourth order time-stepping method,  now symbolized as the ``$2\odot 2=4$" method. The intuitive meaning is that we adopt the second order flow solvers as building blocks and use a two-stage time-stepping to achieve fourth order accurate numerical methods for hyperbolic problems or convection-dominated problems. We make a diagram in the following.
\begin{center}
\fbox{$``4"$: A fourth order scheme}  = \\[3mm]
\fbox{$``2"$: Second order Lax-Wendroff type flow solvers}\\
 + \\
 \fbox{$``2"$: A two-stage time stepping  }
\end{center}
\vspace{0.2cm}

 Careful readers may observe the validity of \eqref{eq:semi} when the above two-stage algorithm applies to the current case, which is why we have to develop the Lax-Wendroff type flow solvers based on hyperbolic balance laws rather than the formal partial differential equations \cite{M-Li-2018}.  There are at least two points that we should concern: {\em (i)  System \eqref{eq:semi} is index-dependent and therefore each equation for fixed $j$ is related to the neighboring  equations; (ii) the continuity of $\mathcal{L}_j$ is crucial when applying the above two-stage algorithm.  }
It is at these points that  \eqref{eq:dy} is substantially different from \eqref{eq:ode}.  Therefore, the regularity of flux is a key factor to guarantee the validity of this algorithm.  Physically speaking, this regularity is natural by recalling the Lagrangian form of fluid dynamical systems \cite{Landau-1987}. Hence in the computation of instantaneous values \eqref{eq:pair-3} or \eqref{eq:pair-2d}, we must be aware of the regularity of the flux that will be further emphasized in the next section about the GRP solver.

 \section{ The GRP solver: A discontinuous and nonlinear  LW flow solver}\label{sec:GRP}

 As is well-known, and also pointed out in the last section, the standard Lax-Wendroff solver results in an algorithm that producing oscillatory solutions if discontinuities are present.
  The GRP solver, the abbreviation of the generalized Riemann problem (GRP) solver,   can be regarded as the discontinuous and nonlinear  version of the Lax-Wendroff  solver. This solver was originally proposed in \cite{Ben-Artzi-84} for compressible fluid flows and related problems. See \cite{Ben-Artzi-Book} for the comprehensive summary of works before $2003$. Later on a direct  Eulerian version of GRP solver was derived in \cite{M-Li-2006} and further extended to general hyperbolic conservation laws \cite{M-Li-2007,Qian-Li-2014}.  The presentation below will follow  the direct Eulerian  GRP.  The application to non-conservative systems is referred to \cite{Balsara-Li-2018}.

 \vspace{0.2cm}

\subsection{1-D  GRP solver}

We first review one-dimensional GRP solver for hyperbolic balance laws
\begin{equation}
\bu_t+\bbf(\bu)_x=\bh(x,\bu),
\label{eq:balance}
\end{equation}
subject to the initial data of form \eqref{data:1D}. An important prototype is   the compressible Euler equations with cross section,
\begin{equation}
\begin{array}{l}
\dfr{\pt (A(x)  \rho)}{\pt t}+\dfr{\pt (A(x)\rho u)}{\pt x} =
0,\\\\
\dfr{\pt (A(x)\rho u)}{\pt t}+\dfr{\pt (A(x)\rho u^2)}{\pt x}+A(x)\dfr{\pt p}{\pt x}
=0,\\\\
\dfr{\pt (A(x)\rho E)}{\pt t}+\dfr{\pt ( A(x)u(\rho E+p))}{\pt x} =0,\\
\end{array}
\label{duct}
\end{equation}
where the variables $\rho$, $u$, $p$ and $E$ are the density,
velocity,  pressure and the total specific energy. The total
specific energy consists of two parts $E=\frac{u^2}2+e$, $e$ is
the internal specific energy. The function $A(x)$ is the area of
the duct. When $A(x)\ev 1$, the system (\ref{duct}) represents the
planar compressible Euler equations.
Let $T$ be the temperature. Then
the entropy $S$ can be defined, as usual, by Gibbs relation of
thermodynamics,
\begin{equation}
TdS = de- \dfr{p}{\rho^2}d\rho. \label{second-law}
\end{equation}
The local sound speed $c$ is defined as
\begin{equation}
c^2 =\dfr{\pt p(\rho, S)}{\pt \rho}.
\end{equation}

We will distinguish  the linear (acoustic) and nonlinear GRP solvers. Both are related. However, as strong waves are involved, the nonlinear GRP solver becomes crucial.  Details can be found first in \cite{Matania-F-duct-1986} and later in \cite{M-Li-2007,Qian-Li-2014} for  general versions.
\vspace{0.2cm}

We denote
\begin{equation}
\begin{array}{l}
\bu_\ell =\bP_-(0-0),\ \ \ \  \ \bu_r  =\bP_+(0+0), \\[3mm]
\bu_\ell' =\bP'_-(0-0), \ \ \ \  \ \bu_r  =\bP_+'(0+0).
\end{array}
\end{equation}
\vspace{0.2cm}

 \n {\bf (I)  Linear GRP solver.}    Linear hyperbolic  equations describe the propagation of waves linearly. For the single linear   case \eqref{eq:adv},
 \begin{equation}
 u_t +a u_x =\al u,
 \end{equation}
with a damping constant $\al$,  we have obviously
 \begin{equation}
 u_0 =\dfr{a+|a|}2 u_\ell +  \dfr{a-|a|}2 u_r ,  \ \ \ \  \left(\dfr{\pt u}{\pt t}\right)_0 =-\dfr{a+|a|}2 u_\ell'- \dfr{a-|a|}2 u_r '+\al u_0.
 \end{equation}
 \vspace{0.2cm}

 For the linear or semi-linear system case
 \begin{equation}
 \bu_t +\bA \bu_x=\bh(\bu,x),
 \end{equation}
 we need to diagonalize the system and pursue  the characteristic decomposition.  Denote by $\la_1,\cdots, \la_m$ the eigenvalues of $\bA$, $\bm\La =\mbox{diag}(\la_1,\cdots, \la_m)$, and  $|\bm\La| =\mbox{diag}(|\la_1|,\cdots, |\la_m|)$.  Then we decompose $\bA$ as $ \bA =\bR \bm\La \bR^{-1}$, and denote again $|\bA| = \bR |\bm \La|\bR^{-1}$. It turns out that the instantaneous values take
 \begin{equation}
 \begin{array}{l}
 \bu_0 = \dfr{\bA+|\bA|}{2}\bu_\ell +  \dfr{\bA-|\bA|}{2} \bu_r , \\[3mm]
    \left(\dfr{\pt \bu}{\pt t}\right)_0 =-\dfr{\bA+|\bA|}2 \bu_\ell'- \dfr{\bA-|\bA|}2 \bu_r '+\bh(u_0,0).
    \end{array}
 \end{equation}
 Therefore, the GRP solver for  the linear system case is substantially the same as that for  the single equation case.
 \vspace{0.2cm}

\n  {\bf (II) Acoustic Approximation}

 For nonlinear cases,  if $\bu_\ell=\bu_r $, but $\|\bu_\ell'-\bu_r '\| \neq 0$, only linear waves emanate from the singularity point $(0,0)$. Then  we can linearize the system \eqref{eq:balance},   at the value $\bu_0 = \bu_\ell=\bu_r $, as
 \begin{equation}
 \bm\th_t + \bA(\bu_0) \bm\th_x =\bh(x, \bu_0), \ \ \ \bm\th =\bu-\bu_0,
 \label{eq:linearized}
 \end{equation}
where $\bA(\bu_0)$ is the Jacobian $\bbf'(\bu_0)$.   Note that
 \begin{equation}
  \left(\dfr{\pt \bu}{\pt t}\right)_0 = \left(\dfr{\pt \bm\th}{\pt t}\right)_0,  \ \ \mbox{ and } \left(\dfr{\pt \bu}{\pt x}\right)_0 = \left(\dfr{\pt \bm\th}{\pt x}\right)_0 .
 \end{equation}
 Immediately we have
 \begin{equation}
  \left(\dfr{\pt \bu}{\pt t}\right)_0=-\dfr{\bA(\bu_0)+|\bA(\bu_0)|}2 \bu_\ell'- \dfr{\bA(\bu_0)-|\bA(\bu_0)|}2 \bu_r '+\bh(u_0,0).
  \label{eq:deri-1}
 \end{equation}

 We can proceed  to obtain any higher temporal derivatives $(\pt^m \bu/\pt t^m)_0$, abbreviated as ADER in \cite{ADER-2002}.  However, in the framework of $2\odot 2=4$,  we are satisfied with the first order temporal variation in \eqref{eq:deri-1}.
 \vspace{0.2cm}

 As $\|\bu_\ell-\bu_r \|\ll 1$, weakly nonlinear waves emanate from $(0,0)$. Then we can carry out the so-called {\em acoustic approximation.}
 To be precise, we can use either exact or approximate Riemann solvers  to obtain the intermediate state $\bu_0$ and linearize the system \eqref{eq:balance}  to be in the form \eqref{eq:linearized} so that the temporal derivative $(\pt \bu/\pt t)_0$ is calculated as in \eqref{eq:deri-1}.
 \vspace{0.2cm}

 For the Euler equations \eqref{duct}, the acoustic approximation is the following.

 \begin{enumerate}
\item[ (i)] As $u_0-c_0>0$ or $u_0+c_0<0$, the acoustic waves moves to one side of the $t$-axis. Then $(\pt \bu/\pt t)_0$ is taken upwind.
\vspace{0.2cm}

\item[ (ii)]  As $u_0-c_0<0<u_0+c_0$, the $t$-axis is located between two acoustic waves. Then we have
\begin{equation}
\begin{array}{ll}
\d \left(\dfr{\pt u}{\pt t}\right)_0=&
\d -\dfr 12
\left[(u_0+c_0)\left(u_\ell'+\dfr{p_\ell'}{\rho_0 c_0}\right )
+(u_0-c_0) \left(u_r '-\dfr{p_r '}{\rho_0
c_0} \right)\right],\\[3mm]
\d  \left(\dfr{\pt p}{\pt t}\right)_0=& \d -\dfr {\rho_0c_0}2
\left[(u_0+c_0)\left(u_\ell'+\dfr{p_\ell'}{\rho_0 c_0}\right )
-(u_0-c_0) \left(u_r '-\dfr{p_r '}{\rho_0
c_0} \right)\right]\\[3mm]
&\d -\dfr{A'(0)}{A(0)} \rho_0 c_0^2 u_0.\\
\end{array}
\end{equation}
The quantity $(\pt \rho/\pt t)_0$ is solved according to the direction of the contact discontinuity,
\begin{equation}
\left(\dfr{\pt \rho}{\pt t}\right)_0= \dfr{1}{c_0^2} \left[
\left(\dfr{\pt p}{\pt t}\right)_0 +u_0 \left( p_\ell'-c_0^2\rho_\ell'
\right) \right] \label{acoustic2}
\end{equation}
if $u_0=u_\ell=u_r >0$; and
\begin{equation}
\left(\dfr{\pt \rho}{\pt t}\right)_0= \dfr{1}{c_0^2} \left[
\left(\dfr{\pt p}{\pt t}\right)_0 +u_0 \left( p_r '-c_0^2\rho_r '
\right) \right] \label{acoustic3}
\end{equation}
if $u_0=u_\ell=u_r <0$.
\end{enumerate}
 \vspace{0.2cm}

 For  ``cheap"  engineering applications, one can use the local Riemann solution $\bu_0$ to linearize the nonlinear system and obtain a ``linear" system so that the above acoustic approximation strategy can be adopted.

 \vspace{0.2cm}

\n  {\bf (III) Nonlinear  GRP solver}

\begin{figure}[!htb]
  \subfigure[Acoustic GRP solver with different grid points ]{
  \includegraphics[width=2.5in]{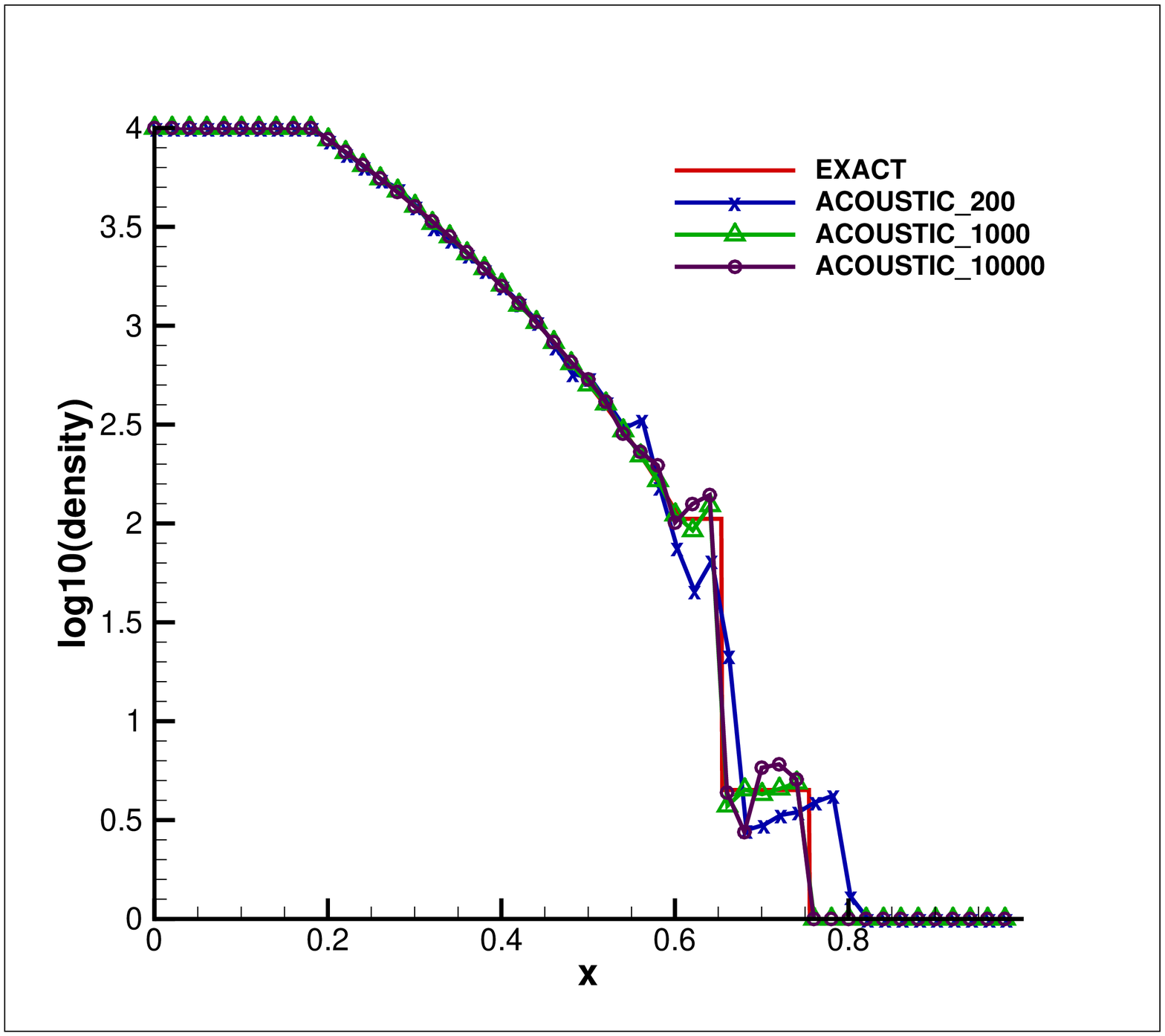}}
\hspace{1cm}
\subfigure[Zoomed solution]{  \includegraphics[width=2.5in]{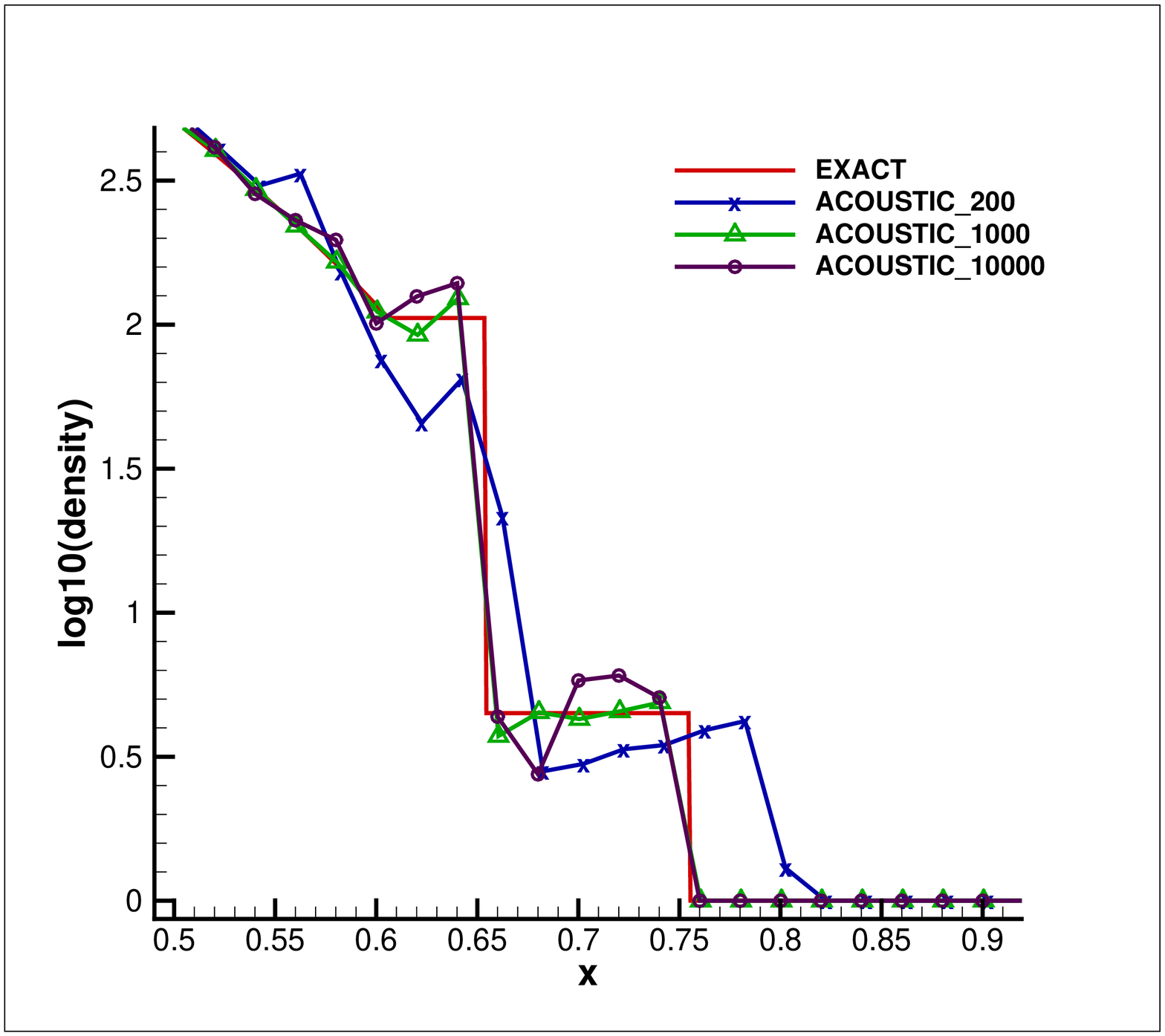} }
  \caption{The numerical solutions computed by the second order acoustic GRP  scheme (with the exact Riemann solver) are compared with the exact solution (only 66 cells are shown). }
  \label{Fig-second-acoustic}
\end{figure}

\begin{figure}[!htb]
  \subfigure[GRP with relatively small number of  grid points ]{
  \includegraphics[width=2.5in]{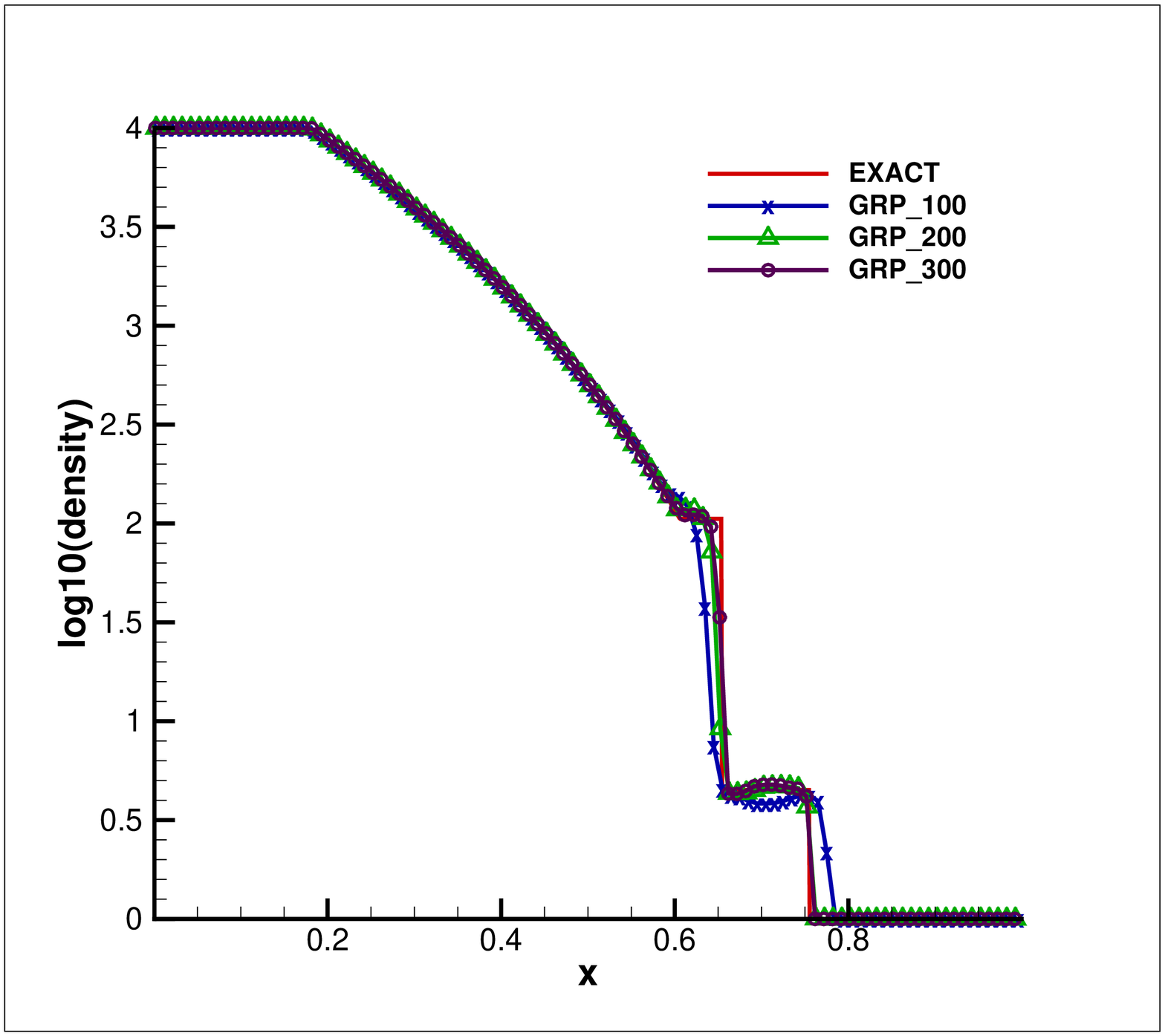}}
\hspace{1cm}
\subfigure[Zoomed solution]{  \includegraphics[width=2.5in]{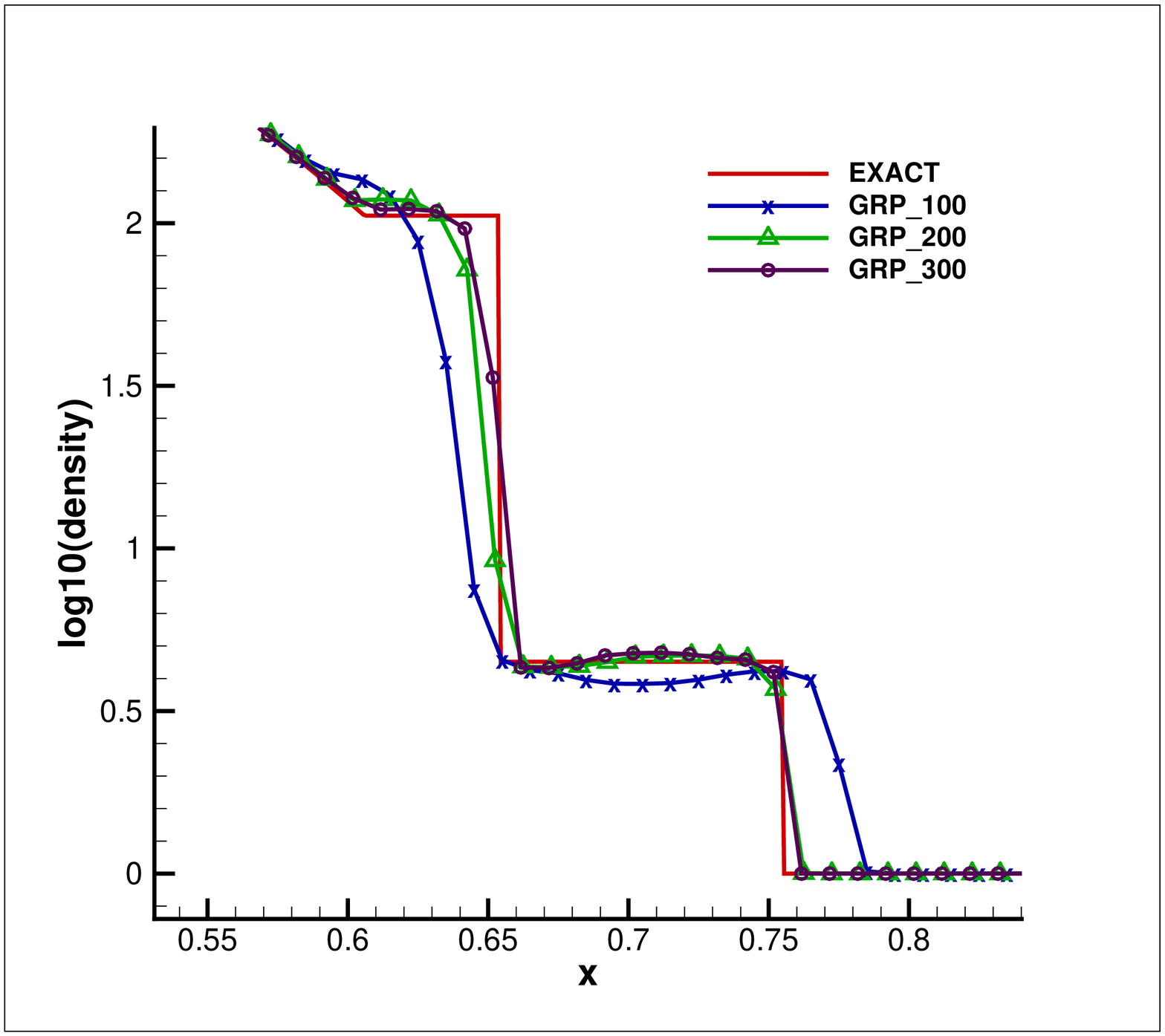} }
\caption[small]{The GRP simulation (only 100 cells are shown)}
\label{Fig-GRP}
\end{figure}

As $\|\bu_\ell-\bu_r \|\gg 1$,  nonlinear waves emanate from the singularity point $(0,0)$. The larger the difference between $\bu_\ell$ and $\bu_r $ is, the stronger the strength of the waves is.
 In Figure \ref{Fig-second-acoustic}  we use the acoustic GRP solver to  simulate the {\em big density ratio problem} and observe  that the numerical solution has large disparity from the exact solution.  In Figure
  \ref{Fig-GRP},  we use the nonlinear GRP solver that will be described below, and see that the numerical solution is improved prominently \cite{Li-Wang-2017}.  Hence it is essential to  develop the nonlinear GRP solver as long as strong waves need resolving.  We just illustrate the nonlinear GRP solver for Euler equations with cross section \eqref{duct}. For general hyperbolic balance laws, readers are  referred to \cite{M-Li-2007,Qian-Li-2014}.

\begin{figure}[!htb]
\centering \subfigure[Wave pattern for the GRP. The initial data
$\bu(x,0)= \bu_\ell+x\bu_\ell'$ for $x<0$ and $\bu(x,0)=\bu_r +\bu_r ' x$ for $x>0$. ]{
\psfrag{0}{$0$}
\psfrag{x}{$x$}\psfrag{t}{$t$}\psfrag{betal}{\small
$\beta=\beta_\ell$}\psfrag{ul}{$\bu_\ell$}\psfrag{al2}{$\bar{\bar\al}$}
\psfrag{shock}{shock}\psfrag{rarefaction}{rarefaction}\psfrag{uminus}{$\bu_-(x,t)$}\psfrag{uplus}{$\bu_+(x,t)$}
\psfrag{u1}{$\bu_{1*}$}\psfrag{u2}{$\bu_{2*}$}\psfrag{al=albar}{$\al=\bar\al$}\psfrag{xbar}{$\bar
\al$ }\psfrag{al=2albar}{$\al=\bar{\bar \al}$}\psfrag{ur}{$\bu_r $}
\psfrag{contact}{contact}\psfrag{betastar}{$\beta=\beta_0$}
\psfrag{ul}{$\bu_\ell$}\psfrag{ur}{$\bu_r $}\psfrag{ustar}{$\bu_0$}
\includegraphics[clip, width=9cm]{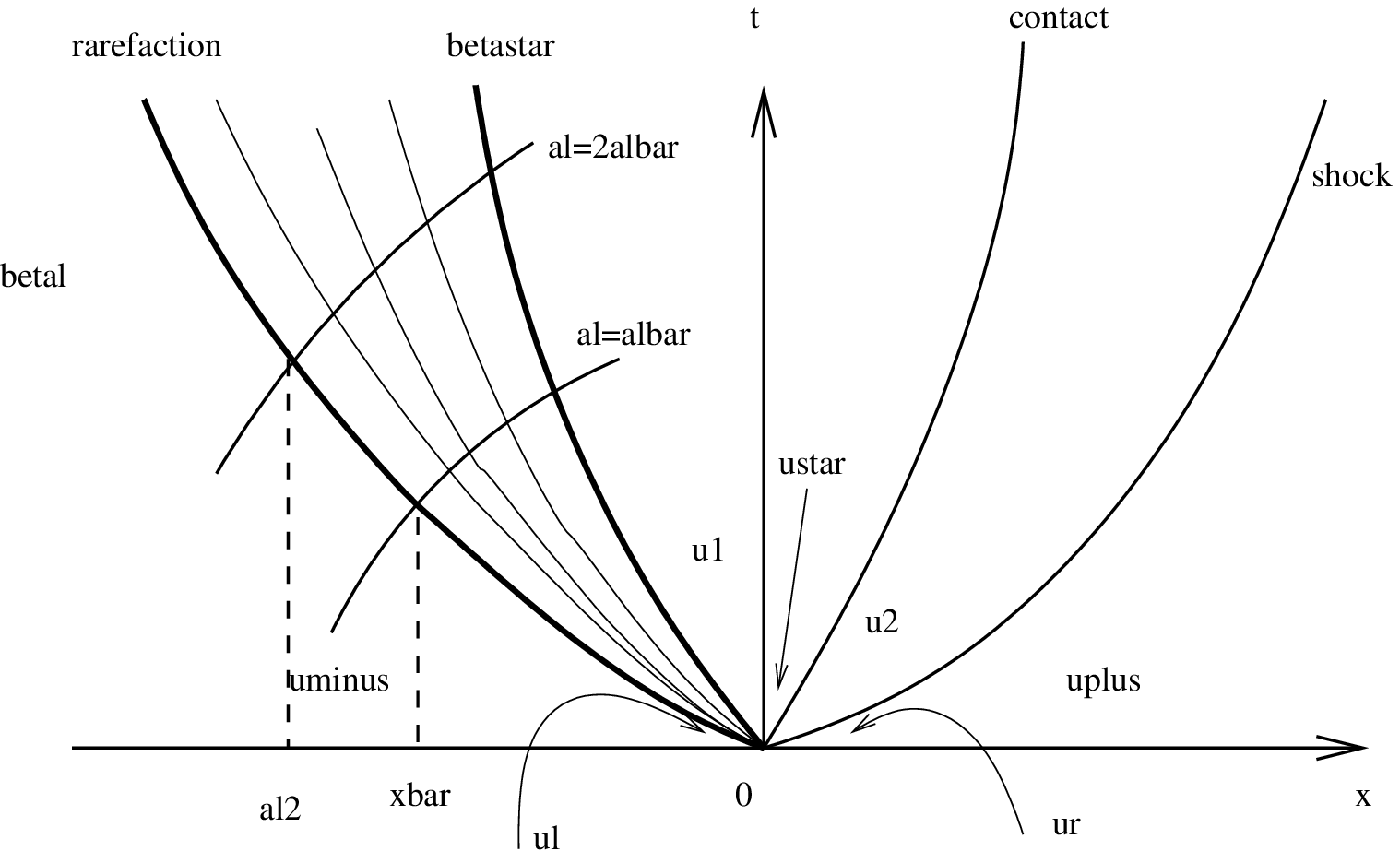} }
\hspace{0.5cm}
 \subfigure[Wave pattern for the associated Riemann problem]{
\psfrag{x}{$x$}\psfrag{t}{$t$}\psfrag{betal}{\small
$\beta=\beta_\ell$}
\psfrag{shock}{shock}\psfrag{rarefaction}{rarefaction}\psfrag{uminus}{$\bu_\ell$}\psfrag{uplus}{$\bu_r $}
\psfrag{u1}{$\bu_{1*}$}\psfrag{u2}{$\bu_{2*}$}\psfrag{al=albar}{$\al=\bar\al$}\psfrag{xbar}{$\bar
\al$ }\psfrag{al=2albar}{$\al=\bar{\bar \al}$}
\psfrag{contact}{contact}\psfrag{betastar}{$\beta=\beta_0$}\psfrag{ustar}{$\bu_0$}
\psfrag{0}{$0$}\psfrag{al2}{$\bar{\bar \al}$}
\includegraphics[width=9cm]{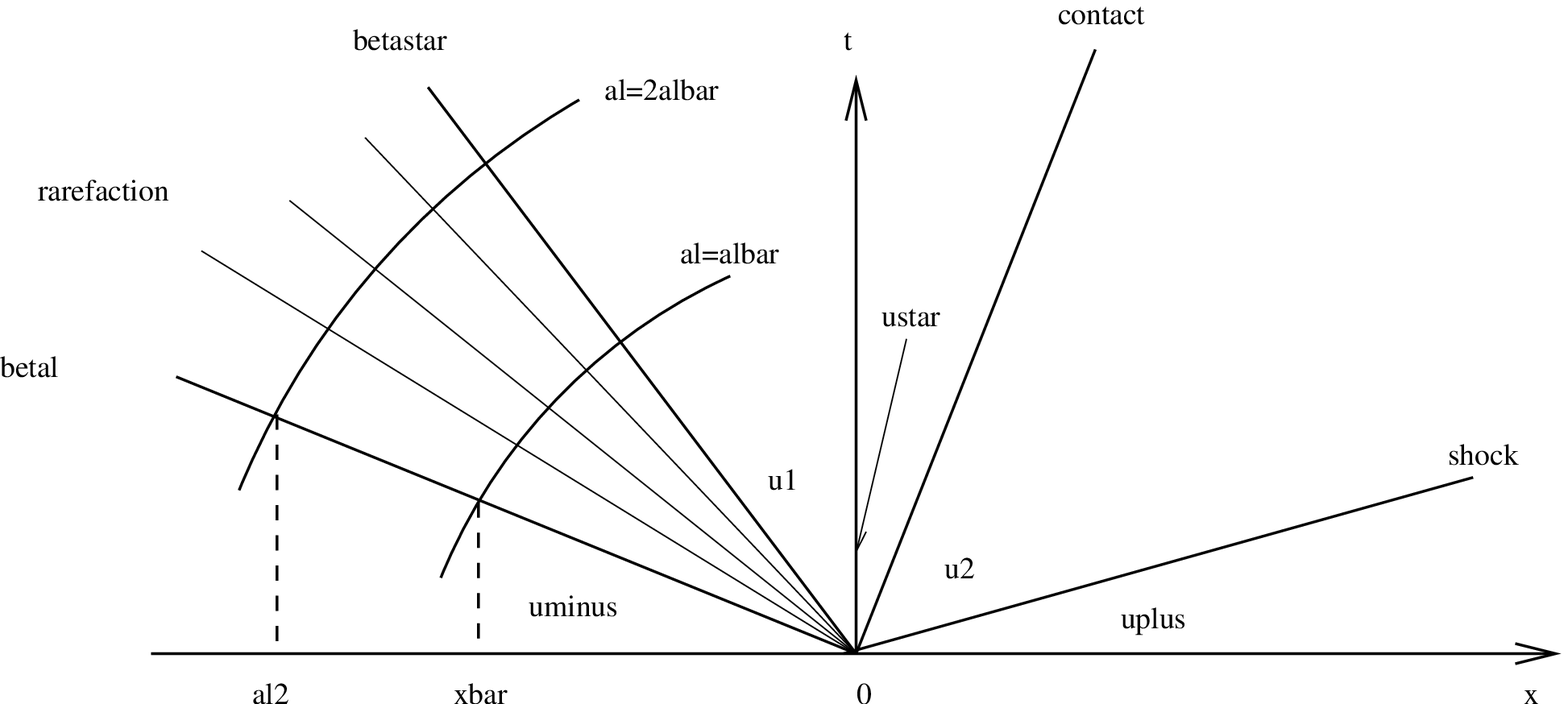}}
\caption[small]{Typical wave pattern for the generalized Riemann problem}
\label{Fig-wave}
\end{figure}

We just need to assume a typical case, as shown in Figure \ref{Fig-wave},  that a rarefaction wave moves to the left, a shock moves to the right and a contact discontinuity lies in the middle.  When the waves move to one side of $t$-axis, i.e., $u_\ell-c_\ell>0$ or $u_r+c_r<0$, it is treated as the acoustic case that $(\pt \bu/\pt t)_0$ is obtained upwind.
\vspace{0.2cm}

We rewrite \eqref{duct} in terms of  $(\rho, u, S)$,
\begin{equation}
\begin{array}{l}
\dfr{\pt \rho}{\pt t} + u\dfr{\pt \rho}{\pt x}+\rho \dfr{\pt
u}{\pt x} =-\dfr{A'(x)}{A(x)}\rho u,\\[3mm]
\dfr{\pt u}{\pt t}+u\dfr{\pt u}{\pt x} +\dfr 1\rho \dfr{\pt p}{\pt
x}=0,\\[3mm]
\dfr{\pt S}{\pt t}+u\dfr{\pt S}{\pt x}=0,
\end{array}
\label{smooth-flow}
\end{equation}
where $p$ is regarded as a function of $\rho$ and $S$. In terms of $\rho$, $u$ and $p$,
the third equation of (\ref{smooth-flow}) can be replaced by,
\begin{equation}
\dfr{\pt p}{\pt t} +u \dfr{\pt p}{\pt x} +\rho c^2 \dfr{\pt u}{\pt
x}= -\dfr{A'(x)}{A(x)} \rho c^2 u.
\end{equation}
In order to resolve strong rarefaction waves, it is particularly essential to introduce the so-called {\em generalized Riemann invariants}, as in \cite{M-Li-2007},
\begin{equation}
\phi= u - \int^\rho \dfr{c(\om,S)}{\om}d\om, \ \ \ \ \psi= u
+\int^\rho \dfr{c(\om, S)}{\om}d\om. \label{def-rie}
\end{equation}
 Together with the entropy variable $S$,   system \eqref{duct}  becomes
\begin{equation}
\left\{
\begin{array}{l}
\dfr{\pt \phi}{\pt t}+(u-c)\dfr{\pt \phi}{\pt x} = B_1,\\[3mm]
\dfr{\pt \psi}{\pt t}+(u+c)\dfr{\pt \psi}{\pt x} =B_2 ,\\[3mm]
\dfr{\pt S}{\pt t} + u\dfr{\pt S}{\pt x}=0.
\end{array}
\right. \label{duct-coup}
\end{equation}
where $B_1=  T\frac{\pt S}{\pt x}+\frac{A'(x)}{A(x)} cu$,
$B_2=T\frac{\pt S}{\pt x}-\frac{A'(x)}{A(x)} cu$.  Here it is easily seen that the variable section $A(x)$ acts on the dynamical behavior of $\phi$ and $\psi$, and thus on that of $(\rho, u, p)$.  The severe change of entropy inevitably  leads to the variation of other physical variables. The GRP solver that we  will derive  tells precisely how the entropy and the cross section affect the dynamics.
\vspace{0.2cm}

The most important ingredient is the application of nonlinear geometric optics for the local expanding of rarefaction waves using local characteristic coordinates $(\al,\beta)$, as shown in Figure \ref{Fig-wave}.  With that,  we can obtain the instantaneous temporal derivatives $\pt S/\pt t$ and $\pt \psi/\pt t$ as  (restricted to polytropic gases),
\begin{equation}
\begin{array}{l}
T\dfr{\pt S}{\pt t}(0,\beta)= - (\beta+c_\ell\th)\th
^{\frac{2\gm}{\gm-1}} T_\ell S_\ell',\
\ \ \ \ \th=c(0,\beta)/c_\ell,\\
\dfr{\pt \psi}{\pt t}(0,\beta) = G_1+\dfr{A'(0)}{2A(0)}G_2,
\end{array}
\label{duct-s-psi}
\end{equation}
where $T_\ell S_\ell'$ is defined through the Gibbs relation,
\begin{equation}
T_\ell S_\ell' =e_\ell' -\dfr{p_\ell}{\rho^2_\ell} \rho_\ell',
\end{equation}
and $G_1$,
$G_2$ are given by,
\begin{equation}
\begin{array}{l}
G_1= -\dfr{\beta+c_\ell\th}{c_\ell} \th^{\frac{\gm+1}{\gm-1}}T_\ell S_\ell' +
 \dfr{\beta+2c_\ell\th}{c_\ell} \th^{\frac{3-\gm}{2(\gm-1)}}\left[\dfr{2\gm}{3\gm-1}T_\ell S_\ell'-c_\ell\psi_\ell'\right],\\\\
G_2=\beta(\beta+c_\ell\th) -(\beta+2c_\ell\th)\left[u_\ell \th^{\frac{3-\gm}{2(\gm-1)}}+ \ol{G_{2}}\right],\\\\
\ol{G_{2}} = \left\{
\begin{array}{ll}
\dfr{-2(\gm+1)c_\ell\th}{(\gm-1)(3\gm-5)}\left(1-\th^{\frac{5-3\gm}{2(\gm-1)}}\right)-
\dfr{(\gm+1)\psi_\ell}{\gm-3}\left[1-\th^{\frac{3-\gm}{2(\gm-1)}}\right],  &\mbox{if } \gm\neq 3,5/3,\\
2c_\ell(\th-1) -\psi_\ell\ln\th, &\mbox{if  } \gm=3,\\
2\left[c_\ell\th\ln\th+\psi_\ell\left(1-\th\right)\right], &\mbox{if }
\gm =5/3.
\end{array}
\right.
\end{array}
\label{duct-coe}
\end{equation}
For general cases, please refer to \cite{Li-Wang-2018}.
\vspace{0.2cm}

Several remarks are in order  about the role of entropy variation and the cross section on the dynamics.

\begin{enumerate}
\item[(i)] The source term reflecting geometric variation always plays an important role in the dynamics of flows. Inherently, the critical  gas indices are clearly exhibited in \eqref{duct-coe}, which cannot be illustrated in other flow solvers.  This is just an evidence that  only  GRP solver can distinguish different gases.

\item[(ii)]  The entropy change rate  is essential  and acts on other flow variables as its initial variation is severe. This tells why the GRP solver is indispensable  when strong waves are simulated. As the involved waves are weak or  $T \frac{\pt S}{\pt x}$  is small, $\frac{\pt S}{\pt t}$ is negligible and many approximations such as linearization are acceptable.

\end{enumerate}

The shock is resolved by tracking its trajectory described by the Rankine-Hugoniot relation,
\begin{equation}
\sg=\dfr{\rho u-\b \rho \b u}{\rho -\b \rho},\ \ \ u= \b u\pm
\Phi(p; \b p, \b \rho), \ \ \ \ \ \rho=\Psi(p; \b p,\b \rho),
\label{iso-RH}
\end{equation}
where $(\rho, u, p)$ and $(\b \rho, \b u,\b p)$ are the states
ahead and behind the shock, respectively, and
\begin{equation}
\Phi(p; \b p,\b \rho)=(p-\b p)\sqrt{\dfr{1-\mu^2}{\b \rho(p+\mu^2
 \b p)}}, \ \ \ \  \ \Psi(p; \b p,\b\rho)=  \b\rho\dfr{p+\mu^2 \b p}{\b p+\mu^2p},
  \ \ \  \  \mu^2 =\dfr{\gm-1}{\gm+1},
\end{equation}
for polytropic gases.
 \vspace{0.2cm}

 For the contact discontinuity $x=x(t)$, we make use of the continuity property of pressure and velocity on both sides of the trajectory,
 \begin{equation}
 u(x(t)-0,t) =u(x(t)+0,t),  \ \ \ \ \ p(x(t)-0,t) =p(x(t)+0,t).
 \end{equation}
 Then the differentiation along the contact discontinuity gives
 \begin{equation}
 \dfr{Du}{Dt} (x(t)-0,t) =\dfr{Du}{Dt} (x(t)+0,t), \ \  \dfr{Dp}{Dt} (x(t)-0,t) =\dfr{Dp}{Dt} (x(t)+0,t),
 \end{equation}
 where $D/Dt=\pt/\pt t+u\pt/\pt x$ is the material derivative. This relation bridges the rarefaction wave and the shock, just like that for the Riemann solver.  We just remind that the density and entropy undergo jump across this contact discontinuity.
 \vspace{0.2cm}

 Thus we come to the nonlinear GRP solver that are distinguished as  nonsonic and sonic cases.

\begin{prop} (Non-sonic case.)  Assume a typical wave configuration for the generalized
Riemann problem for (\ref{duct}) as shown in Figure \ref{Fig-wave} that the $t$-axis is located in the intermediate region. Then  $(\pt u/\pt t)_0$ and
$(\pt p/\pt t)_0$ satisfies  the following pair of linear equations,
\begin{equation}
\begin{array}{l}
a_\ell\left(\dfr{\pt u}{\pt t}\right)_0 +b_\ell\left(\dfr{\pt p}{\pt t}\right)_0= d_\ell, \\
a_r \left(\dfr{\pt u}{\pt t}\right)_0 +b_r \left(\dfr{\pt p}{\pt t}\right)_0= d_r ,
\end{array}
\end{equation}
where $a_\ell$, $b_\ell$, $d_\ell$ and $a_r $, $b_r $, $d_r $ are specified
below. Also the computation of $(\pt \rho/\pt t)_0$ are  computed
by the following two cases. \vspace{0.2cm}

(i) If $u_0>0$, the contact discontinuity moves to the right and separates two states $(\rho_{0\ell}, u_0, p_0)$, $(\rho_{0r}, u_0, p_0)$. The coefficients $a_\ell$, $b_\ell$ and $d_\ell$ are given as,
\begin{equation}
(a_\ell, b_\ell, d_\ell)= (\td a_\ell, \td b_\ell, \td d_\ell).
\end{equation}
The coefficients $a_r $, $b_r $ and $d_r $  are given by
\begin{equation}
\begin{array}{l}
a_r = \dfr{c_{0\ell}^2}{c_{0\ell}^2-u_0^2}\left[\td a_r  (1-\dfr{\rho_{0\ell}u_0^2}{\rho_{0r}c_{0r}^2})
+\td b_r  (\rho_{0r}-\rho_{0\ell})u_0\right],\\\\
b_r =\dfr{1}{c_{0\ell}^2-u_0^2}\left[\td a_r  (-\dfr{1}{\rho_{0\ell}}+\dfr{c_{0\ell}^2}{\rho_{0r}c_{0r}^2})u_0
-\td b_r  (-\dfr{\rho_{0r}}{\rho_{0\ell}} u_0^2+c_{0\ell}^2)\right],\\\\
d_r =\td d_r  +\dfr{A'(0)}{A(0)} \dfr{u_0^3}{c_{0\ell}^2-u_0^2}\left[\td a_r  (1-\dfr{\rho_{0\ell}c_{0\ell}^2}{\rho_{0r}c_{0r}^2})
+\td b_r (\rho_{0r}-\rho_{0\ell})c_{0\ell}^2\right].
\end{array}
\end{equation}
The value $(\pt \rho/\pt t)_0$ is computed from the rarefaction side,
\begin{equation}
\left(\dfr{\pt \rho}{\pt t}\right)_0 =
\dfr{1}{c_{0\ell}^2}\left[\left(\dfr{\pt p}{\pt t}\right)_0+(\gm-1)
\rho_{0\ell}
u_0\left(\dfr{c_{0\ell}}{c_\ell}\right)^{\frac{1+\mu^2}{\mu^2}}T_\ell S_\ell'\right],
\label{rho-rate}
\end{equation}
by using the state equation $p=p(\rho,S)$.
\vspace{0.2cm}

(ii) If $u_0<0$, the contact discontinuity moves to the left. The coefficients $a_r $, $b_r $ and $d_r $ are given
as,
\begin{equation}
(a_r , b_r , d_r )= (\td a_r , \td b_r , \td d_r ).
\end{equation}
While the coefficients $a_\ell$, $b_\ell$ and $d_\ell$  are computed by
\begin{equation}
\begin{array}{l}
a_\ell= \dfr{c_{0r}^2}{c_{0r}^2-u_0^2}\left[\td a_\ell (1-\dfr{\rho_{0r}u_0^2}{\rho_{0\ell}c_{0\ell}^2})
+\td b_\ell (\rho_{0\ell}-\rho_{0r})u_0\right],\\\\
b_\ell=\dfr{1}{c_{0r}^2-u_0^2}\left[\td a_\ell (-\dfr{1}{\rho_{0r}}+\dfr{c_{0r}^2}{\rho_{0\ell}c_{0\ell}^2})u_0
-\td b_\ell (-\dfr{\rho_{0\ell}}{\rho_{0r}} u_0^2+c_{0r}^2)\right],\\\\
d_\ell=\td d_\ell +\dfr{A'(0)}{A(0)} \dfr{u_0^3}{c_{0r}^2-u_0^2}\left[\td a_\ell (1-\dfr{\rho_{0r}c_{0r}^2}{\rho_{0\ell}c_{0\ell}^2})
+\td b_\ell(\rho_{0\ell}-\rho_{0r})c_{0r}^2\right].
\end{array}
\end{equation}
The value $(\pt \rho/\pt t)_0$ is computed from the shock side,
\begin{equation}
g_{\rho}^R \left(\dfr{\pt \rho}{\pt t}\right)_0 +g_p^R
\left(\dfr{Dp}{Dt}\right)_0+g_u^R\left(\dfr{Du}{Dt}\right)_0 = u_0
\cdot f_r , \label{rho-deri-shock}
\end{equation}
where  $g_{\rho}^R$, $g_p^R$, $g_u^R$ and
$f_r $ are constant, explicitly given in the following,
\begin{equation}
\begin{array}{c}
g_{\rho}^R = u_0- \sg_0,   \ \  g_p^R= \dfr{\sg_0}{c_{0r}^2
}-u_0 H_1,  \ \ \ g_u^R= \rho_{2*}(\sg_0-u_0)\cdot u_0 \cdot H_1,\\\\
\begin{array}{rl}
f_r =& (\sg_0-u_r )\cdot H_2 \cdot p_r '+(\sg_0-u_r )\cdot H_3
\cdot\rho_r ' -\rho_r  \cdot \left(H_2 c_r ^2+H_3\right)\cdot
u_r '\\\\
& -\dfr{A'(0)}{A(0)}\cdot \left(H_2 c_r ^2+H_3\right)\rho_r u_r ,
\end{array}\\\\
\end{array}
\end{equation}
and $H_1$, $H_2$ and $H_3$ are expressed by
\begin{equation}
\begin{array}{c}
 H_1=\dfr{ \rho_r (1-\mu^4) p_r }{( p_r +\mu^2 p_0)^2}, \ \ \ H_2=\dfr{
\rho_r (\mu^4-1)p_0}{( p_r +\mu^2 p_0)^2},\ \ \ H_3=\dfr{p_0+\mu^2  p_r }{
p_r +\mu^2 p_0}.
\end{array}
\end{equation}
\vspace{0.2cm}

The other coefficients $(\td a_\ell, \td b_\ell, \td d_\ell)$ and  $(\td a_r, \td b_r, \td d_r)$  are
\begin{equation}
\begin{array}{l}
\td a_\ell(0,\beta)= 1, \\[3mm]  \td b_\ell(0,\beta) =
\dfr{1}{\rho(0,\beta) c(0,\beta)},\\[3mm]
\td d_\ell(\beta)=\dfr{\beta+2\th c_\ell}{c_\ell}\cdot
\th^{\frac{3-\gm}{2(\gm-1)}}\left(\dfr{2\gm}{3\gm-1}T_\ell S_\ell'-c_\ell\psi_\ell'\right)+\dfr{A'(0)}{2A(0)}
G_2, \label{duct-dl}
\end{array}
\end{equation}
and
\begin{equation}
\begin{array}{l}
\td a_r =1-\dfr{\sg_0 u_0}{u_0^2-c_{0r}^2}-\dfr{\sg_0 \rho_{0r}c_{0r}^2}{u_0^2-c_{0r}^2}\cdot \Phi_1,\\[3mm]
\td b_r = \dfr{1}{\rho_{0r}}\dfr{\sg_0 }{u_0^2-c_{0r}^2}-\left(1-\dfr{\sg_0 u_0}{u_0^2-c_{0r}^2}\right)\Phi_1,\\[3mm]
\td d_r =L_p^R \cdot  p_r '+L_u^R \cdot u_r '+L_{\rho}^R \cdot \rho_r ' +\dfr{A'(0)}{A(0)} j_r ,
\end{array}
\end{equation}
where we use notations
\begin{equation}
\begin{array}{l}
L_p^R= -\dfr{1}{\rho_r }+(\sg_0-u_r )\cdot \Phi_2 , \\
L_u^R=
\sg_0-u_r  -\rho_r \cdot c_r ^2 \cdot \Phi_2-\rho_r  \cdot \Phi_3,\\
L_{\rho}^R= (\sg_0-u_r )\cdot \Phi_3, \\
j_r  = -(\Phi_2 c_{0r}^2+\Phi_3)\rho_r  u_r  +
\left(1+\Phi_1\rho_{0r}u_0\right)\dfr{\sg_0 c_{0r}^2 u_0}{u_0^2-c_{0r}^2}; \\\\
\sg_0=\dfr{\rho_{0r}u_0-\rho_r u_r }{\rho_{0r}-\rho_r },\\
\Phi_1 =\dfr 12 \sqrt{\dfr{1-\mu^2}{ \rho_r (p_0+\mu^2  p_r )}}
\cdot \dfr{p_0+(1+2\mu^2) p_r }{p_0+\mu^2  p_r }, \\
\Phi_2=-\dfr 12 \sqrt{\dfr{1-\mu^2}{ \rho_r (p_0+\mu^2  p_r )}}
\cdot
\dfr{(2+\mu^2)p_0+\mu^2  p_r }{p_0+\mu^2  p_r } , \\
\Phi_3 =- \dfr{p_0- p_r }{2 \rho_r } \sqrt{\dfr{1-\mu^2}{
\rho_r (p_0+\mu^2 p_r )}}.
\end{array}
\label{right-coef}
\end{equation}
\end{prop}

\vspace{0.2cm}

\begin{prop}(Sonic case). \label{prop-duct-sonic} Assume that the $t-$axis is located
inside the rarefaction wave associated with $u-c$.
Then we have
\begin{equation}
\begin{array}{l}
\left(\dfr{\pt u}{\pt t}\right)_0 =\dfr 12\left[ \td d_\ell+\th^{\frac{2\gm}{\gm-1}} T_\ell S_\ell'+\dfr{A'(0)}{A(0)} u_0^2\right],\\[3mm]
\left(\dfr{\pt p}{\pt t}\right)_0 =\dfr {\rho_0 c_0}2
\left[ \td d_\ell-\th^{\frac{2\gm}{\gm-1}} T_\ell S_\ell'-\dfr{A'(0)}{A(0)} u_0^2\right],
\end{array}
\label{sonic-case}
\end{equation}
where $\td d_\ell$ is given in (\ref{duct-dl}), with $\th=c_0/c_\ell$,
and $(u_0,\rho_0, c_0)$ is the limiting value of $(u,\rho,c)$
along the $t$-axis so that $u_0-c_0=0$. Then density change rate is given as in \eqref{rho-rate}.
\end{prop}

The above formulae look complicated, but seem irreplaceable.  We can go to \cite{M-Li-2007} for technical derivation of them.

 \subsection{Temporal-spatial coupling and thermodynamical effect}

 Thermodynamics distinguishes compressible fluid flows from incompressible ones, and the Mach number can be regarded as a parameter of  the compressibility. The entropy dissipation is a necessary condition guaranteeing the stability of numerical schemes. Let's consider  the  compressible Euler equations with uniform cross section ($A(x)\ev constant$ in \eqref{duct}).
The entropy inequality says
\begin{equation}
(\rho S)_t + (\rho u S)_x \geq  0, \ \ \mbox{ in } \mathcal{D}'.
\end{equation}
However, it is a well-known open problem whether this inequality is satisfied at discrete level, particularly for high order accurate schemes. There are two origins of discrete errors: the data projection and the flux approximation.
In a general setting of finite volume framework,  given the initial data $\bu_n(x) \in \mathcal{P}_k$ at $t=t_n$, we have to find the solution $\bu_{n+1}(x)$ at next time level $t=t_{n+1}$, satisfying
\begin{equation}
\begin{array}{rl}
\d \int_{x_{j-\frac 12}}^{x_{j+\frac 12}} \rho S(\bu_{n+1}(x))dx & \d \geq \int_{x_{j-\frac 12}}^{x_{j+\frac 12}} \rho S(\bu_n(x)) dx\\[3mm]
&\d  -\De t\left[ (\rho u S)_{j+\frac 12}^{approx}  - ( \rho u S)_{j-\frac 12}^{approx} \right]+ Tol(\De x, \De t),  \\[4mm]
&\d  .
\end{array}
\label{eq:entropy}
\end{equation}
where $(\rho uS)_{j+\frac 12}^{approx}$ is the numerical entropy flux, and $Tol(\De x,\De t)$ is the entropy production that has the maximum tolerance of order three, $Tol(\De x,\De t)=\mathcal{O}(\De t^3+\De x^3)$.
 We comment on how to achieve this inequality at the discrete level in the following.

\begin{enumerate}
\item[(i)]  The persistence space $\mathcal{P}_k$ often consists of piecewise polynomials of degree $k$. Given the initial data $\bu_n(x)\in \mathcal{P}_k$, we solve \eqref{duct} and obtain the (analytic) entropy solution $\bu(x,t)$ for $t_n<t<t_{n+1}$, satisfying
\begin{equation}
\begin{array}{rl}
&\d \int_{x_{j-\frac 12}}^{x_{j+\frac 12}}\rho S(\bu(x,t_{n+1}))dx  \geq  \d \int_{x_{j-\frac 12}}^{x_{j+\frac 12}}\rho S(\bu_n(x))dx\\[4mm]
-&\d  \left[\int_{t_n}^{t_{n+1}}\rho u S(x_{j+\frac 12},t)dt-\int_{t_n}^{t_{n+1}}\rho u S(x_{j-\frac 12},t)dt\right].
\end{array}
\end{equation}

\item[(ii)] The projection of $\bu(x,t)$ onto $\mathcal{P}_k$ (reconstruction procedure) is required to satisfy
\begin{equation}
\int_{x_{j-\frac 12}}^{x_{j+\frac 12}}\rho S(\bu_{n+1} (x))dx \geq \int_{x_{j-\frac 12}}^{x_{j+\frac 12}}\rho S(\bu(x,t_{n+1}))dx+\mathcal{O}(\De x^3).
\label{entropy-P}
\end{equation}
This is an extremely difficult  step. For scalar conservation laws, there was a nice discussion on MUSCL-type linear reconstruction \cite{Bouchut-1996}.  In general, the Jensen inequality tells that
\begin{equation}
\begin{array}{l}
\d \rho S(\bar\bu_j^{n+1}) \geq \dfr{1}{\De x} \int_{x_{j-\frac 12}}^{x_{j+\frac 12}}\rho S(x,t_{n+1})dx, \\[3mm]
\d  \bar\bu_j^{n+1} =\dfr{1}{\De x}  \int_{x_{j-\frac 12}}^{x_{j+\frac 12}}\bu(x,t_{n+1})dx.
\end{array}
\end{equation}
Hence, there is still some room to make \eqref{entropy-P} hold, which remains an open problem.
\vspace{0.2cm}

\n{\em Open problem on data projection: Find an optimal reconstruction strategy so that \eqref{entropy-P} holds.}

\item[(iii)]  Assume that \eqref{entropy-P} holds for certain data projection.  As shown above, we approximate the interface value in the following way
\begin{equation}
\bu(x_{j+\frac 12},t) = \bu_{j+\frac 12}^n +\left( \dfr{\pt \bu}{\pt t}\right)_{j+\frac 12}^n(t-t_n) +\mathcal{O}((t-t_n)^2), \ \ \ t_n<t<t_{n+1}.
\end{equation}
 In particular, we make  use of  the entropy information. It turns out that
\begin{equation}
(\rho u S)_{j+\frac 12}^{approx} : = (\rho u S)(x_{j+\frac 12},t_{n+\frac 12} )= \dfr{1}{\De t} \int_{t_n}^{t_{n+1}} \rho u S(x_{j+\frac 12},t)dt+\mathcal{O}(\De t^2).
\label{eq:en-flux}
\end{equation}
Summarizing all together yields \eqref{eq:entropy}.
\end{enumerate}
\vspace{0.2cm}

It is observed   that the precise calculation of entropy in \eqref{duct-s-psi} is a direct  way to achieve \eqref{eq:en-flux}.  Other ways may at most lead to
\begin{equation}
(\rho u S)_{j+\frac 12}^{approx} : = (\rho u S)(x_{j+\frac 12},t_{n+\frac 12} )= \dfr{1}{\De t} \int_{t_n}^{t_{n+1}} \rho u S(x_{j+\frac 12},t)dt+\mathcal{O}(\|\bu\|^2).
\label{eq:en-flux1}
\end{equation}
The error of $\mathcal{O}(\|\bu\|^2)$ is not tolerated unless for scalar cases or smooth flows, since this type of errors violate the entropy inequality in the limit.
\vspace{0.2cm}

It is no doubt that the achievement of \eqref{eq:en-flux} is the outcome of the direct use of the entropy equation in \eqref{smooth-flow}, which is actually the Lax-Wendroff procedure, a temporal-spatial coupling procedure.

 \subsection{M-D GRP solver and transversal effects} When computing multidimensional (M-D)  problems, M-D GRP solver is necessary to reflect the transversal effect, which is impossible  using the exact or approximate normal Riemann solvers. We restrict to two-dimensional hyperbolic balance laws
 \begin{equation}
 \bu_t +\bbf(\bu)_x +\bg(\bu)_y=\bh(x,y,\bu),
 \label{law:2D}
 \end{equation}
 where $\bbf$ and $\bg$ are flux functions. 3-D GRP solver is straightforward.   The initial data takes the form
  \begin{equation}
 \bu_0(x,y) = \left\{
 \begin{array}{ll}
 \bu_-(x,y),  \ \ & x<0,\\[3mm]
 \bu_+(x,y), & x>0,
 \end{array}
 \right.
 \label{data:2D-h}
 \end{equation}
 where $\bu_\pm(x,y)$ are two polynomials of degree $k$. The $x-$direction is the normal and the $y$-direction is the transversal.  A particular case is
 \begin{equation}
\bu_0(x,y) =\left\{
\begin{array}{ll}
 \bu_- +\nb\bu_{0,-}\cdot \bx, \ \ \ x<0,\\[3mm]
\bu_+ +\nb\bu_{0,+}\cdot \bx, \ \ \ x>0.
\end{array}
\right.
\label{data:2D-2}
\end{equation}
 The M-D GRP solver can be classified as {\em M-D linear GRP solver, acoustic GRP solver, nonlinear normal GRP solver with transversal perturbation, and genuinely nonlinear M-D GRP solver}.
 \vspace{0.2cm}

 \n {\bf (I) M-D linear GRP solver.}   We consider the linear case
 \begin{equation}
 \bu_t +\bA \bu_x +\bB\bu_y=0,
 \label{law:2D-linear}
 \end{equation}
 where $\bA$ and $\bB$ are constant matrices, and both of them have their respective real eigenvalues and the complete sets of eigenvectors.    We first assume that \eqref{law:2D-linear} is subject to the initial data \eqref{data:2D-2}. Note that $\nb \bu$ satisfies the same form of  \eqref{law:2D-linear},
\begin{equation}
(\nb \bu)_t +\bA(\nb\bu)_x+\bB(\nb\bu)_y=0,
\label{g-2}
\end{equation}
but subject to the initial  data
 \begin{equation}
\nb\bu_0(x,y)=\left\{\begin{array}{ll}
\nb\bu_{0,-}, & x<0,\\[3mm]
\nb\bu_{0,+}, & x>0.
\end{array}
\right.
\label{g-1}
\end{equation}
This boils down to the standard Riemann problem for $\nb \bu$. Therefore, the gradient $\nb \bu$ has an explicit expression,
\begin{equation}
\nb\bu(x,y,t) =   \sum_{\ell=1}^{L(x,t)} \nb v_{0,+}^\ell \br^\ell + \sum_{\ell=L(x,t)+1}^m \nb v_{0,-}^\ell \br^\ell,
\label{g-3}
\end{equation}
where  the notations $L(x,t)$ is the maximum value of $\ell$ such that $x-\la^\ell t>0$, $\la^\ell$ is the eigenvalue of $\bA$ and $\br^\ell$ is the associated eigenvector, $\bv$ is so defined that
\begin{equation}
\bu =\sum_{\ell=1}^m v^\ell \br^\ell,  \ \ \ \bv=(v^1, \cdots, v^m).
\end{equation}
See \cite{LeVeque-book}.
In particular, we have
\begin{equation}
\nb\bu(0,y,t) =   \sum_{\ell: \la^\ell<0} \nb v_{0,+}^\ell \br^\ell + \sum_{\ell:\la^\ell>0} \nb v_{0,-}^\ell \br^\ell.
\label{g-4}
\end{equation}
Then we immediately obtain
 \begin{equation}
 \begin{array}{rl}
\left(\dfr{\pt \bu}{\pt t}\right)_0& \d : = \lim_{t\rw 0^+} \dfr{\pt \bu}{\pt t}(0,y, t) \\[3mm]
& \d  =  -\bA \lim_{t\rw 0^+} \dfr{\pt \bu}{\pt x}(0,y, t)-\bB\lim_{t\rw 0^+} \dfr{\pt \bu}{\pt y}(0,y, t)\\[3mm]
& =-(\bA,\bB)\cdot \nb\bu(0,y,0^+).
\end{array}
 \label{g-5}
 \end{equation}
 \vspace{0.2cm}

As far as the more general initial data \eqref{data:2D-h} is concerned,
 the solution $\bu$ consists of  piecewise  polynomials of the same degree as the initial data since \eqref{law:2D-linear} is a linear system.  Here we are satisfied with the second order GRP solver and solve \eqref{law:2D-linear} at any point $(0,y_0, 0)$ to obtain
 \begin{equation}
\left(\dfr{\pt \bu}{\pt t}\right)_0 =-(\bA,\bB)\cdot \nb\bu(0,y_0,0^+),
 \label{g-6}
 \end{equation}
 where $\nb\bu(0,y_0, 0^+)$ are calculated as the same procedure as above.
 \vspace{0.2cm}

 \n{\bf (II) M-D acoustic GRP solver.}

 For nonlinear cases \eqref{law:2D}, if the initial data \eqref{data:2D-h} is continuous but discontinuous in its derivatives, acoustic waves emanate from the interface $x=0$, just as one-dimensional case. For this case, we might as well take the initial data \eqref{data:2D-2} and  assume that $\bu_-=\bu_+$ but $\|\nb \bu_{0,-}-\nb\bu_{0,+}\|\neq 0$.  Then we linearize \eqref{law:2D} around the state $\bu_0=\bu_-=\bu_+$ to obtain
 \begin{equation}
 \bm \th_t +\bbf'(\bu_0)\bm \th_x+\bg'(\bu_0)\bm \th_y=\bh(x,y,\bu_0), \ \ \ \bm \th =\bu-\bu_0.
 \end{equation}
 Then we can exactly follow the linear case to calculate $(\pt \bu/\pt t)_0$.
 \vspace{0.2cm}

 The acoustic approximation applies for the case $\|\bu_--\bu_+\|\ll 1$.  We linearize \eqref{law:2D} around the intermediate state $\bu_0$ resulting from the associated Riemann problem.  Then the linear GRP solver applies for this case.

 \vspace{0.2cm}

 \n{\bf (III) M-D genuinely nonlinear GRP solver with transversal description.}
 \vspace{0.2cm}

 As $\|\bu_--\bu_+\|\gg 1$, we have to deal with genuinely nonlinear GRP.   Thinking of the initial value problem for \eqref{law:2D} subject to the initial data \eqref{data:2D-2}, the solution is the envelope of Riemann solution along $x=0$ locally at $t=0$.  Hence along $x=0$, the associated Riemann solution is known. For example, at two points $(0, y_1,0)$ and $(0,y_2,0)$, we solve the Riemann problem locally for the normal conservation law, respectively,
 \begin{equation}
 \bu_t+\bbf(\bu)_x=0,
 \end{equation}
 and obtain the local intermediate values $\bu(0,y_1, 0^+)$ and $\bu(0,y_2,0^+)$. Then for any point $(0,y,0)$, $y_1<y<y_2$, we can approximate $\frac{\pt \bu}{\pt y}(0,y,0^+)$. Particularly, we approximate
 \begin{equation}
 \bu_y(0,y_0, 0^+) =\dfr{\bu(0,y_1,0^+)-\bu(0,y_2,0^+)}{y_1-y_2}+\mathcal{O}((y_1-y_2)^2).
 \end{equation}
 Then we regard the transversal term $\bg(\bu)_y$ and $\bh(x,y,\bu)$ as a perturbation locally, and solve the following problem,
 \begin{equation}
 \bu_t +\bbf(\bu)_x =-\bg(\bu)_y +\bh(x,y,\bu)=: \bd(x,y,\bu,\bu_y).
 \end{equation}
 This boils down to the one-dimensional planar problem locally, for which the GRP solver was proposed in \cite{M-Li-2007}.
Detailed and complete M-D GRP solver is proposed in \cite{Qi-PhD-2017}.

 \subsection{Transversal effects for genuinely M-D schemes}

 For multidimensional (M-D) problems,  the balance law can be always written in the form,
 \begin{equation}
 \dfr{d}{dt} \int_{\Om}\bu(\bx,t)d\bx =-\int_{\pt\Om} \bbf(\bu)\cdot \bn d L,
 \end{equation}
 where $\Om$ is the control volume, $\pt \Om$ is the boundary and $\bn$ is the outer unit normal.  We ignore the external force just for the clarity of presentation.
 \vspace{0.2cm}

 Thanks to the Galilean invariance, we always assume that $(1,0)$  (the direction $x$-axis) is the normal direction, and $(0,1)$ (the direction of $y$-axis) is the transversal direction.   The standard Riemann solver just reflects  the normal effect. In contrast, the LW procedure can describe the transversal effect precisely. Consider a linear advection problem
\begin{equation}
u_t+a u_x+bu_y=0.
\end{equation}
Then we use the temporal-spatial coupling property to obtain
\begin{equation}
\left(\dfr{\pt u}{\pt t} \right)_{\pt \Om} =-a \left(\dfr{\pt u}{\pt x} \right)_{\pt \Om}-b\left(\dfr{\pt u}{\pt y} \right)_{\pt \Om},
\end{equation}
where $(\frac{\pt u}{\pt x})_{\pt\Om}$ and  $(\frac{\pt u}{\pt y})_{\pt\Om}$ can be obtained by solving the associated Riemann problem.  Also as remarked in Section 2 for the linear wave system,  the transversal effect is substantial even though the convergence rate is formally the same.

 \section{The kinetic LW flow solver}\label{sec:GKS}

The fluid dynamics can be described in various viewpoints, such as the kinetic description. The governing equation is the Boltzmann-type equation
\begin{equation}
f_t +\bxi\cdot\nb_\bx f =\dfr{1}{\ep} B(f, f),
\label{Boltz}
\end{equation}
where $f=f(t,\bx, \bxi)$ is the density distribution, $\bxi$ is the velocity of molecules (particle),  and $B(f, f)$ is the collision term, $\ep$ is the Knusner number.
Ideally, for a given initial distribution, we solve \eqref{Boltz} to obtain the solution $f(t,x_{j+\frac 12}, \bxi)$  and define the numerical flux
\begin{equation}
\bbf_{j+\frac 12}(t_n; t_{n+1}) =\dfr{1}{\De t} \int_{t_n}^{t_{n+1}}\int_{\Re^3} \bxi \psi f(t,x_{j+\frac 12}, \xi)d\bxi dt,
\end{equation}
 where $\psi=(1,\bxi,\bxi^2)^\top$ is the invariant, and the average of macroscopic variables is
 \begin{equation}
 \bu_j^{n+1}=\dfr{1}{\De x} \int_{x_{j-\frac 12}}^{x_{j+\frac 12}}\int_{\Re^3} \psi f(t_{n+1},x,\bxi)d\bxi dx.
 \end{equation}
 In general, it is difficult to solve the equation \eqref{Boltz} analytically.  To understand the relation of macroscopic equations (Euler and Navier-Stokes equations) and the kinetic equation, we first take the {\em so-called railroad method} in \cite{Ohwada-2002} as an example to illustrate how to devise kinetic schemes.

\subsection{Railroad method}
Consider the linear advection equation
\begin{equation}
u_t + a u_x=0, \ \ \  u(x,0) =u_0(x).
\label{eq:ad}
\end{equation}
Introduce a distribution function $f(t,x,\xi)$ and define the macroscopic variable $u(x,t)$ as a moment of $f$,
\begin{equation}
u(x,t) =\int_{-\iy}^{\iy} f(t,x,\xi)d\xi.
\label{data:ad}
\end{equation}
If we choose $f$ to take the form,
\begin{equation}
f(t,x,\xi)= \dfr{u(x,t)}{\sq{\pi}}\exp[-(\xi-a)^2],
\end{equation}
then $f(t,x,\xi)$ satisfies
\begin{equation}
f_t +\xi f_x = Q[f]:= \dfr{(\xi-a)}{\sq{\pi}}\exp[-(\xi-a)^2]\dfr{\pt u}{\pt x},
\label{rail-f}
\end{equation}
subject to initial data
\begin{equation}
f(0,x,\xi) =\dfr{u_0(x)}{\sq{\pi}}\exp[-(\xi-a)^2].
\label{rail-d}
\end{equation}
The initial value problem \eqref{eq:ad} and the problem \eqref{data:ad}-\eqref{rail-d} are equivalent:  If one solution is known, then the other is defined.
We write the solution of  \eqref{rail-f}-\eqref{rail-d} as
\begin{equation}
f(t,x,\xi) = f(0, x-\De t\xi,\xi) +\int_0^{\De t} Q[f](\tau,x-(\De t-\tau)\xi),\tau,\xi)d\tau.
\label{sol:f}
\end{equation}
Then the solution $u(x,\De t)$ is given as
\begin{equation}
u(x,t) = \int_{-\iy}^{\iy}f(0, x-\De t\xi,\xi) d\xi +\int_{-\iy}^\iy\int_0^{\De t} Q[f](\tau,x-(\De t-\tau)\xi),\tau,\xi)d\tau d\xi.
\end{equation}
Note that the LW approach for \eqref{sol:f} yields
\begin{equation}
f(t,x,\xi) = f(0,x,\xi) + \De t \xi \frac{\pt}{\pt x} f(0,x,\xi) + \De t Q[f](0, x,\xi)+ \mathcal{O}(\De t^2).
\end{equation}
Therefore  we have
\begin{equation}
u(x,t) = \int_{-\iy}^\iy (f(0,x,\xi) + \De t \xi \dfr{\pt}{\pt x}  f(0,x,\xi) + \De t Q[f](0, x,\xi) )d\xi +\mathcal{O}(\De t^2),
\end{equation}
which yields a second order approximation to the exact  solution $u(x,t)$.\vspace{0.2cm}

The numerical solution is
\begin{equation}
u_j(\De t) =u_j(0) -\dfr{\De t}{\De x} [F_{j+\frac 12}-F_{j-\frac 12}],
\label{scheme:ad}
\end{equation}
where
$$
\begin{array}{l}
u_j(t)=\dfr{1}{\De x} \int_{x_{j-\frac 12}}^{x_{j+\frac 12}} \int_{-\iy}^\iy f(t,x,\xi)d\xi dx, \\[4mm]
\d F_{j+\frac 12}=\int_0^{\De t} \int_{-\iy}^\iy \xi f(t,x_{j+\frac 12},\xi)d\xi dt.
\end{array}
$$
\vspace{0.2cm}

We assume the initial data for \eqref{eq:ad} is piecewise smooth with possible discontinuity at $x=x_{j+\frac 12}$.  Correspondingly, the initial data \eqref{rail-d} for \eqref{rail-f} consists of two parts.  It turns out that the numerical flux $F_{j+\frac 12}$ in \eqref{scheme:ad} becomes
\begin{equation}
F_{j+\frac 12} =F_{j+\frac 12}^+ +F_{j+\frac 12}^-,
\end{equation}
where $F_{j+\frac 12}^\pm$ consist of three parts, respectively,
\begin{equation}
\begin{array}{l}
\d F_{j+\frac 12}^\pm = \int_0^{\De t} \int_{\pm \xi>0}  \xi f(\tau,x_{j+\frac 12},\xi)d\xi d\tau = \De t G_{j+\frac 12}^\pm -\dfr{\De t^2}2 (H_{j+\frac 12}^\pm -K^\pm_{j+\frac 12}), \\[4mm]
\d G_{j+\frac 12}^\pm =\int_{\pm \xi>0} \xi f(0,x_{j+\frac 12}\mp 0,\xi)d\xi,\\[3mm]
\d H_{j+\frac 12}^\pm =\int_{\pm \xi>0} \xi^2\frac{\pt}{\pt x} f(0,x_{j+\frac 12}\mp 0,\xi)d\xi,\\[3mm]
\d K_{j+\frac 12}^\pm = \int_{\pm\xi>0} \xi Q[f](0,x_{j+\frac 12}\mp 0,\xi)d\xi.
\end{array}
\end{equation}
As the solution is smooth, the scheme \eqref{scheme:ad} becomes the LW approach immediately.  See \cite{Ohwada-2002} for details.

\subsection{The LW type solver for gas kinetic schemes}

Let's now work on a simplified model, the Bhatnagar-Gross-Krook (BGK) model \cite{BGK},
\begin{equation}
f_t + \bm\xi f_\bx =\dfr{g-f}\ep,
\label{eq:BGK}
\end{equation}
where $\ep$ is the collision time, and $g$ is the equilibrium state, approached by $f$ as $\ep$ goes to zero,
\begin{equation}
g= \dfr{\rho}{m} \left(\dfr{m}{2\pi k T}\right)^{\frac 32} e^{-(m/2kT)\bm\xi^2},
\end{equation}
where $m$ and $k$ are constant, $\rho$ and $T$ are density and temperature, respectively.  Indeed, all  macroscopic variables $\rho$, $\bu$ and $E$ are defined as
\begin{equation}
(\rho, \bu, E)(\bx,t) =\int_{\Re^3}  (1,\bm\xi, \bm\xi^2) f(t,\bx,\bm\xi)d\bm\xi.
\end{equation}
The validity of BGK-model is clearly explained in \cite{BGK}.
\vspace{0.2cm}

Starting with \eqref{eq:BGK},  Xu and his collaborators successfully  developed gas kinetic scheme (GKS) solver \cite{Xu-1993,Xu-1994,Xu-2001,Xu-2010,Xu-2015}.  A key ingredient is that the explicit  solution formula for \eqref{eq:BGK} is used for the numerical flux approximation,
\begin{equation}
f(t, x_{j+\frac 12}, \xi) =  \dfr{1}{\ep} \int_{0}^{\De t} g(\tau, x', \xi)e^{-(t-t')/\ep} dt'+e^{-t/\ep} f_0(x_{j+\frac 12}-\xi t, \xi),
\label{gks-solution}
\end{equation}
subject to the initial data $f_0(x,\xi)$, where $x'=x_{j+\frac 12}-\xi(t-t')$. Here just the case of one-dimension is described. The full information contained in \eqref{gks-solution} provides  ``exact"  expression of flux across the interface $x=x_{j+\frac 12}$, which is of course consistent with the LW flow solver. We can go to \cite{Xu-2015} for comprehensive description of the GKS solver.
For the gas-kinetic scheme, the gas evolution is a relaxation process from kinetic to hydrodynamic scale through the exponential function,
and the corresponding flux is a complicated function of time.
\vspace{0.2cm}

In order to obtain the time derivatives of the
flux function at $t_n$ and $t_*=t_n + \Delta t/2$ with the correct physics,
the flux function should be approximated as a linear function of time within a time interval.
Let's first introduce the following notation,
\begin{align*}
\mathbb{F}_{i+1/2}(W^n,\delta)
=\int_{t_n}^{t_n+\delta}F_{i+1/2}(W^n,t)dt&=\int_{t_n}^{t_n+\delta}\int
\xi f,tmx_{i+1/2},\xi)d\xi dt.
\end{align*}
In the time interval $[t_n, t_n+\Delta t]$, the flux is expanded as
the following linear form
\begin{align}\label{expansion}
F_{i+1/2}(W^n,t)=F_{i+1/2}^n+ \partial_t F_{j+1/2}^n(t-t_n).
\end{align}
The coefficients $F_{j+1/2}^n$ and $\partial_tF_{j+1/2}^n$ can be
determined as follows,
\begin{align*}
F_{i+1/2}(W^n,t_n)\Delta t&+\frac{1}{2}\partial_t
F_{i+1/2}(W^n,t_n)\Delta t^2 =\mathbb{F}_{i+1/2}(W^n,\Delta t) , \\
\frac{1}{2}F_{i+1/2}(W^n,t_n)\Delta t&+\frac{1}{8}\partial_t
F_{i+1/2}(W^n,t_n)\Delta t^2 =\mathbb{F}_{i+1/2}(W^n,\Delta t/2).
\end{align*}
By solving the linear system, we have
\begin{align}\label{second}
F_{i+1/2}(W^n,t_n)&=(4\mathbb{F}_{i+1/2}(W^n,\Delta t/2)-\mathbb{F}_{i+1/2}(W^n,\Delta t))/\Delta t,\nonumber\\
\partial_t F_{i+1/2}(W^n,t_n)&=4(\mathbb{F}_{i+1/2}(W^n,\Delta t)-2\mathbb{F}_{i+1/2}(W^n,\Delta t/2))/\Delta
t^2.
\end{align}
Similarly, $\displaystyle F_{i+1/2}(W^*,t_*), \partial_t
F_{i+1/2}(W^*,t_*)$ for the intermediate state can be constructed.  For the two-dimensional
computation, the corresponding fluxes in the $y$-direction can be
obtained as well. Readers are referred to \cite{Xu-Li-2017}.
\vspace{0.2cm}

There are huge numbers of references about  kinetic solvers, which are beyond the scope of the current paper. We stop to discuss further.

 \section{Compact reconstruction using the Hermite interpolation}
\label{sec:inflow}

 The compactness is a key factor in the design  of high order schemes,  determining  the dissipation of the schemes near discontinuities and the numerical treatment of boundary conditions.  With the increase of time-stepping, the width of computational stencils is inevitably expanded for multi-stage methods.  Hence it is very important to construct the data in a compact way.

 Unlike  WENO  using the Lagrangian interpolation, we adopt the Hermite-type interpolation using both the average values of physical (conservative or primitive) variables, and the approximate gradient of the solution. Going back to the original GRP, we construct the data over the computation cell $(x_{j-\frac 12},x_{j+\frac 12})$ as
 \begin{equation}
 u_n(x) =  \bar u_j^n+\sg_j^n (x-x_j),  \ \ \ x\in (x_{j-\frac 12},x_{j+\frac 12}),
 \label{eq:data-linear}
 \end{equation}
 where the gradient is chosen through the procedure,
 \begin{equation}
 \sg_j^n =\dfr{1}{\De x} \minmod(\al (u_j^n-u_{j-1}^n), u_{j+\frac 12}^{n,-} -  u_{j-\frac 12}^{n,-},\al (u_{j+1}^n-u_{j}^n)),  \ \ \al\in [0,2).
 \label{eq:minmod}
 \end{equation}
 Usually, $\al$ is chosen to be as large as possible. As $\al\in (1,2)$, $u_n(x)$ behaves as  sawtooth and  implies that $\sg_j^n$ take mostly
 \begin{equation}
 \sg_j^n =\dfr{1}{\De x} (u_{j+\frac 12}^{n,-} -  u_{j-\frac 12}^{n,-}) =\dfr{1}{\De x} \int_{x_{j-\frac 12}}^{x_{j+\frac 12}} \dfr{\pt u}{\pt x} (x,t_n-0)dx\approx \left(\overline{\dfr{\pt u}{\pt x}}\right)_j^n.
 \label{eq:grad}
 \end{equation}
 This is a natural approximation to the gradient.  The boundary value $u_{j+\frac 12}^{n,-}$, as a strong solution along the cell interface $x=x_{j+\frac 12}$, is calculated from the history,
 \begin{equation}
 u_{j+\frac 12}^{n,-} = u_{j+\frac 12}^{n-1} +\De t\left(\frac{\pt u}{\pt t}\right)_{j+\frac 12}^{n-1},
 \label{eq:inte}
 \end{equation}
 where $u_{j+\frac 12}^{n-1}$ and $\left(\frac{\pt u}{\pt t}\right)_{j+\frac 12}^{n-1}$ are obtained already from the GRP solver, and no extra efforts need making.
 \vspace{0.2cm}

 Some remarks  on \eqref{eq:minmod} are made here, and they can be applied for later high order data interpolations.
 \begin{enumerate}
 \item[(i)] Compared to the classical limiter algorithm, \eqref{eq:minmod} takes \eqref{eq:grad} in smooth regions, and limit the gradient near discontinuities in order to suppress possible oscillations.  This sawtooth-type reconstruction can produce sharper profiles of discontinuities.

 \item[(ii)] The piecewise linear data \eqref{eq:data-linear} is  the embryonic form of Hermite polynomials for high order schemes. Since all values are already given using the GRP solver, no extra effort is made on the calculation of the gradient, unlike in DG methods or other Hermite interpolation \cite{Qiu-Shu-2004}.  If one might argue the freedom of solution elements,  he could regard the current treatment as the Lagrangian interpolation using five points values.

 \item[(iii)]  For the data \eqref{eq:data-linear}, five points values are used for the data reconstruction. Essentially, these values are defined  in three computational cells, rather than in five cell, so that computational stencils are almost half saved.   This is one of key factors achieving compactness.

 \end{enumerate}

 Now we extend this to the two-stage fourth order method, by reviewing the result in \cite{Du-Li-2018}.
   Given the average $\bar u_j$ and the derivative $\De u_j$ of the function $u(x)$  over the cell $I_j$,
\begin{equation}\label{eq:0-th}
  \bar{u}_j = \dfr 1h \int_{I_j}u(x,t)dx,\ \ \ \   \De u_j = \dfr 1h \int_{I_j}\dfr{\pt u}{\pt x}(x,t)dx,
\end{equation}
we want to construct a polynomial  $p(x)$  such that $u_{j+\frac 12,-}$ is its left limiting values at $x=x_{j+\frac 12}$.
Choose three stencils
\begin{equation}\label{eq:stencils}
  S^{(-1)} = I_{j-1}\cup I_j, \ \ S^{(0)} = I_{j-1}\cup I_j \cup I_{j+1}, \ \ S^{(1)} = I_j \cup I_{j+1}.
\end{equation}
On stencil $ S^{(0)} $, $ \bar{u}_{j-1} $, $ \bar{u}_{j} $ and $ \bar u_{j+1} $ are used to construct a polynomial $ p^{(0)} $  for the interpolation. Hence at $ x_{j+\frac 12} $, we have
\begin{equation}\label{eq:s-0}
  u^{(0)}_{j+\frac 12, -} := p^{(0)}(x_{j+\frac 12}) = -\dfr 16 \bar{u}_{j-1} + \dfr{5}{6} \bar{u}_{j} + \dfr {1}{3} \bar u_{j+1}.
\end{equation}
Similarly, $ p^{(-1)} $ and $ p^{(1)} $ are constructed by using $ \bar{u}_j $, $ \bar{u}_{j-1} $, $ \De u_{j-1} $ on $ S^{(-1)} $ and  by using $ \bar{u}_{j} $, $ \bar{u}_{j+1} $, $ \De u_{j+1} $ on $ S^{(1)} $, respectively,
\begin{equation}\label{eq:s-n1-p1}
\begin{array}{l}
  u^{(-1)}_{j+\frac 12, -} := p^{(-1)}(x_{j+\frac 12}) = -\dfr 76 \bar{u}_{j-1} + \dfr{13}{6} \bar{u}_{j} - \dfr {2h}{3} \De u_{j-1},\\[3mm]
  u^{(1)}_{j+\frac 12, -} := p^{(1)}(x_{j+\frac 12}) = \dfr 16 \bar{u}_{j} + \dfr{5}{6} \bar{u}_{j+1} - \dfr {h}{3} \De u_{j+1}.
\end{array}
\end{equation}
If the solution is smooth on the large stencil $ I_{-1}\cup I_0 \cup I_1 $, we have
\begin{equation}\label{eq:s-total}
  \tilde u_{j+\frac 12,-} = \dfr{1}{120}(-23\bar{ u}_{j-1}+76\bar{ u}_{j}+67\bar{u}_{j+1}-9h\Delta u_{j-1}-21h\Delta u_{j+1}).
\end{equation}
Thus the linear weights of the three stencils are
\begin{equation}\label{eq:linear-weight}
  \gm^{(-1)} = \dfr {9}{80}, \ \ \gm^{(0)} = \dfr{29}{80}, \ \ \gm^{(1)} = \dfr{21}{40},
\end{equation}
which ensure
\begin{equation*}
  \tilde\bu_{j+\frac 12,-} = \d\sum_{r=-1}^1\gm^{(r)}\bu^{(r)}_{j+\frac 12, -}.
\end{equation*}
The smoothness indicators are defined by
\begin{equation}\label{eq:SI-def}
  \beta^{(r)} = \d\sum_{l=1}^2 \int_{I_j} h^{2l-1}\left(\dfr{d^l}{dx^l}p^{(r)}(x)\right)^2dx, \ \ r=-1,0,1,
\end{equation}
in the same way as in the WENO reconstructions where $ p^{(r)}(x) $ is the interpolation polynomial on stencil $ S^{(r)} $. Their explicit expressions are
\begin{equation}\label{eq:SI-res}
\begin{array}{l}
  \beta^{(-1)} = (-2\bar{u}_{j-1}+2\bar{u}_{j}-h\De u_{j-1})^2 + \dfr{13}{3}(-\bar{ u}_{j-1}+\bar{u}_{j}-h\Delta u_{j-1})^2,\\[3mm]
  \beta^{(0)} = \dfr 14 (-\bar{u}_{j-1}+\bar{u}_{j+1})^2 + \dfr{13}{12}(-\bar{u}_{j-1}+2\bar{u}_{j}-\bar{u}_{j+1})^2,\\[3mm]
  \beta^{(1)} = (2\bar{u}_{j+1}-2\bar{u}_{j}-h\Delta u_{j+1})^2 + \dfr{13}{3}(\bar{u}_{j+1}-\bar{u}_{j}-h\Delta u_{j+1})^2.
\end{array}
\end{equation}
 Then we compute the nonlinear weights in the same way as the WENO-Z method does
\begin{equation}\label{eq:nonlinear-weight}
  \omega^\text{z}_r = \dfr{\al^\text{z}_r}{\sum_{l}\al_l}, \ \ \ \al^\text{z}_r = \gm^{(r)}(1+\dfr{\tau^\text{z}}{\beta^{(r)}+\varepsilon}), \ \ \ r = -1,0,1,
\end{equation}
where $ \tau^\text{z} = |\beta^{(1)}-\beta^{(-1)}| $ and $ \varepsilon $ is a small parameter in order to avoid a zero denominator. Finally we have
\begin{equation}\label{eq:final}
  u_{j+\frac 12, -} = \d\sum_{r=-1}^1 \omega^\text{z}_r u^{(r)}_{j+\frac 12, -}.
\end{equation}
The right interface value $ \bu_{j-\frac 12, +} $ can be reconstructed in a similar way by mirroring the above procedure with respect to $ x_j = \dfr 12(x_{j-\frac 12} + x_{j+\frac 12}) $.

Since the GRP solver has to use the spatial derivative $(\pt u/\pt x)_{j+\frac 12,\pm}$,  we approximate them using the interpolation,
\begin{equation}\label{eq:ux-lagrange}
  \Big(\dfr{\pt u}{\pt x}\Big)_{j+\frac 12,\pm} := \dfr{1}{12h}\left(\bar u_{j-1}-15\bar u_{j}+15\bar u_{j+1}-\bar u_{j+2}\right).
\end{equation}
 It is observed in this interpolation  does not need the WENO-type stencil selection procedure.

 We define this procedure as {\em HWENO} , terming  a Hermite type interpolation using the WENO interpolation strategy. {\em  GRP4-HWENO5} refers to the  two-stage fourth order scheme based on this Hermite type interpolation using  the GRP solver.

For two-dimensional cases, we can develop the similar approach over rectangular meshes.  See \cite{Du-Li-2018,Ji-Xu-2018}.  As far as unstructured meshes are concerned,  there still remains space to explore.

\begin{figure}[!htp]
\centering
\includegraphics[width=\textwidth]{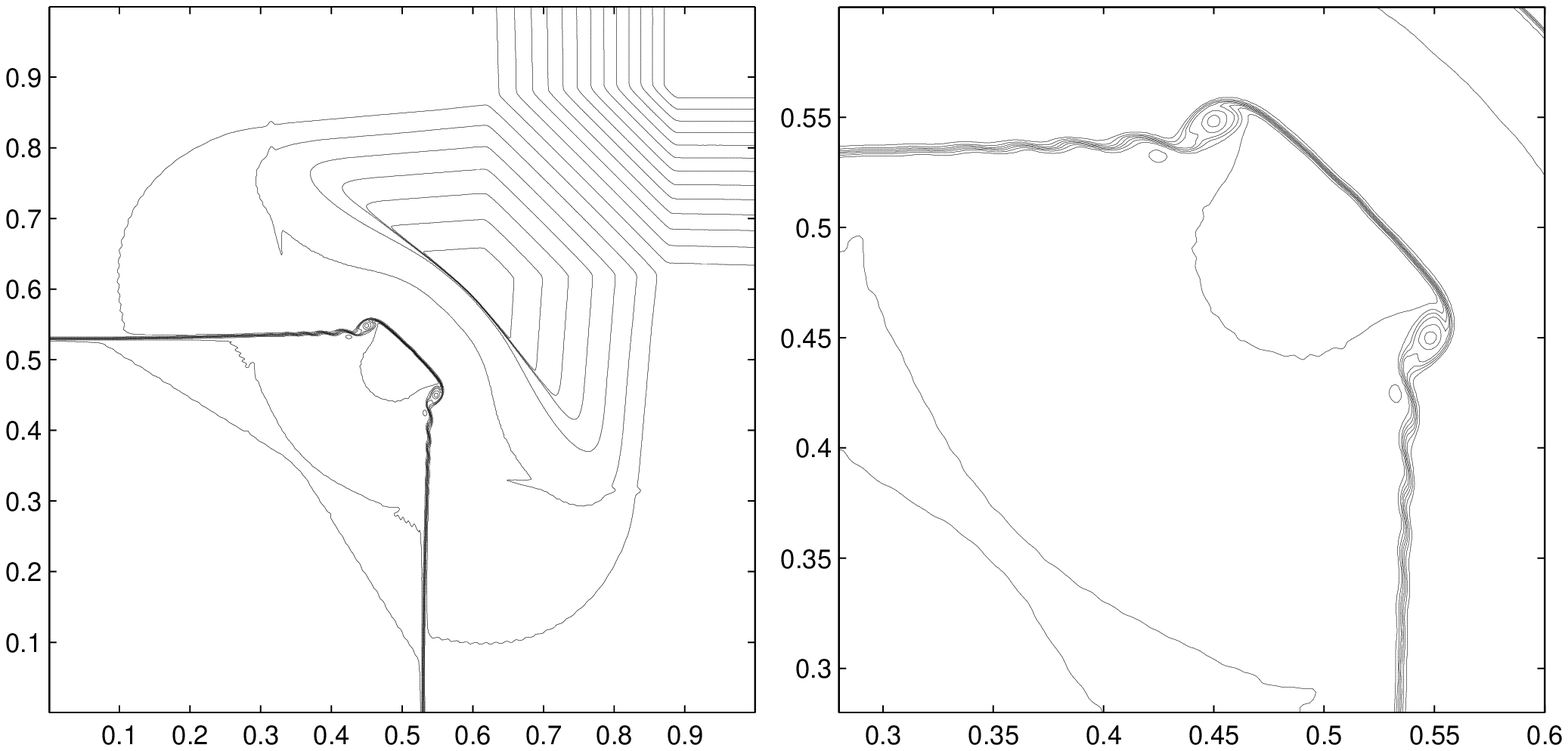}
\includegraphics[width=\textwidth]{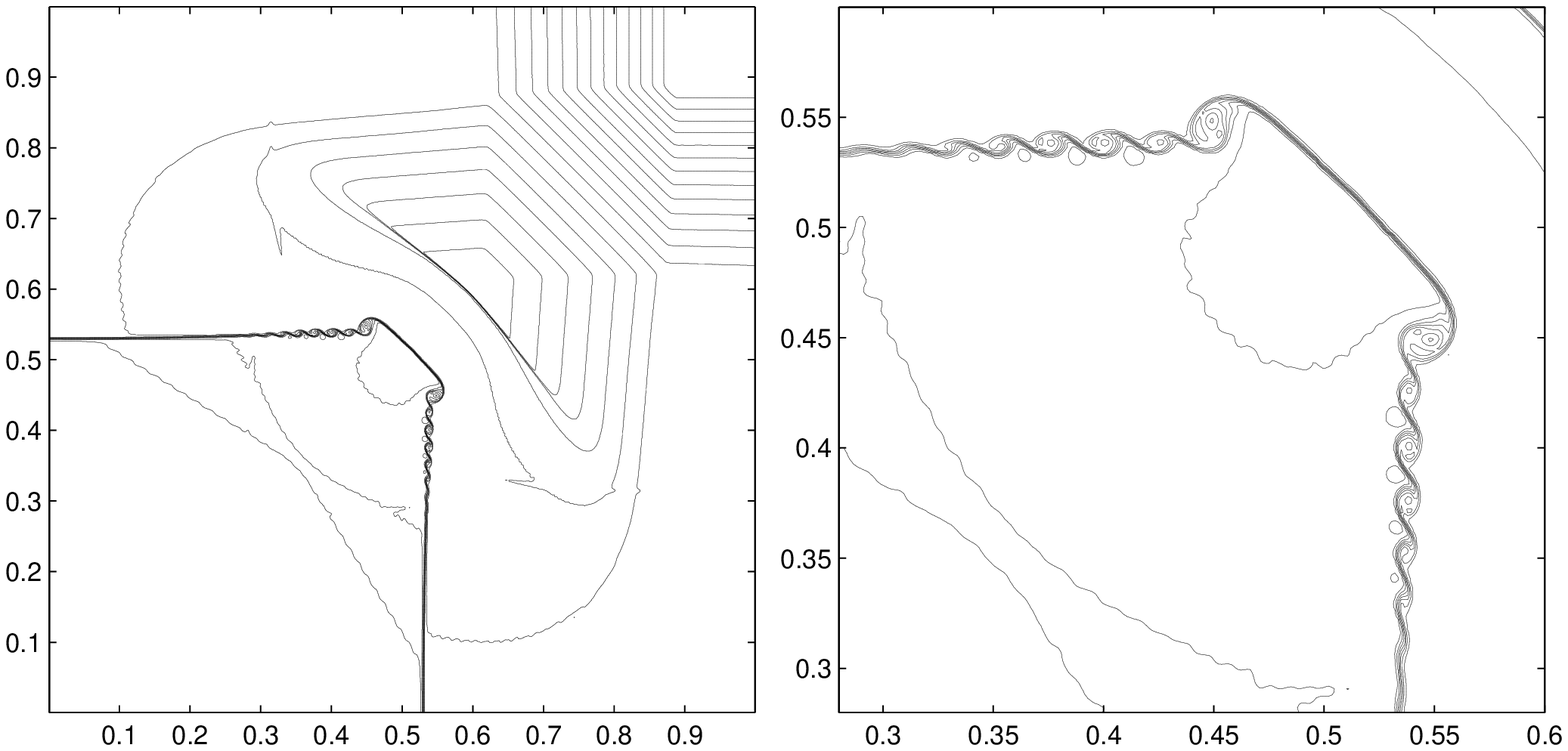}
\caption{The density contours of three 2-D Riemann problem in Example 5 computed with the schemes GRP4-WENO5 (upper) and GRP4-HWENO5 (lower), respectively. $ 700\times700 $ cells are used.}
\label{fig:Riemann-a}
\end{figure}

Here we give an example to demonstrate how important the compactness is.
We provide an example of two-dimensional Riemann problem taken from \cite{Han-Li-2011} involving the interactions of vortex sheets with rarefaction waves. The computation is implemented over the domain $[0,1]\times [0,1]$.
\begin{equation}
(\rho,u,v,p)(x,y,0) =\left\{
\begin{array}{ll}
 (1, 0.1, 0.1, 1), & 0.5<x<1,0.5<y<1,\\
 (0.5197, -0.6259, 0.1, 0.4), & 0<x<0.5, 0.5<y<1,\\
 (0.8, 0.1, 0.1, 0.4),&0<x<0.5,0<y<0.5,\\
 (0.5197, 0.1, -0.6259, 0.4), &0.5<x<1,0<y<0.5.
\end{array}
\right.
\end{equation}
The output time is $0.3$.
The contours of the density and their local enlargements are shown in Figures \ref{fig:Riemann-a}. We can see that the scheme with the Hermite type reconstruction can resolve more small structures along the vortex sheet.

 \section{High Order Boundary Conditions}\label{sec-BC}

 Approximation to boundary conditions may be one of the most challenging  issues in CFD. On one hand, mathematical modelings  of  fluid flows near physical boundaries are diverse.  On the other hand, highly nonlinear behaviors and
complex boundaries make the approximation  notoriously involved.
\vspace{0.2cm}

\begin{figure}[!htb]
\centering
 \includegraphics[scale=0.5]{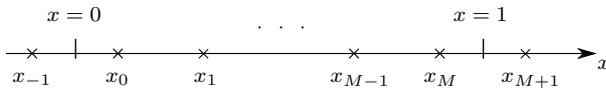}
\small
    \put(-206,15){$ x=0 $}
    \put(-219,-8){$ x_{-1} $}
    \put(-184,-8){$ x_{0} $}
    \put(-151,-8){$ x_{1} $}
    \put(-54,15){$ x=1 $}
    \put(-100,-8){$ x_{M-1} $}
    \put(-65,-8){$ x_{M} $}
    \put(-36,-8){$ x_{M+1} $}
    \put(0,-3){$x$}
 \caption{The computational domain $(0,1)$. Set $ x_0=h/2 $ and $ x_M=1-h/2 $. Then $ x_{-1}=-h/2 $ and $ x_{M+1}=1+h/2 $ are ghost points.}
 \label{fig:kreiss}
\end{figure}

 We briefly illustrate their idea in the finite difference framework by considering  the initial boundary value problem (IBVP) for a scalar conservation law
\begin{equation}
\label{eq:IBVP}
\left\{
\begin{array}{ll}
  \dfr{\pt u}{\pt t} + \dfr{\pt f(u)}{\pt x}=0,\ \  \ &  x \in (0,1),~ t > 0,\\[3mm]
  u(x,0) = u_0(x), & x \in (0,1),\\[3mm]
  u(0,t) = g(t), & t > 0.
\end{array}
\right.
\end{equation}
Assume that $ f^\prime(u)>0 $ for all $ u\in\Re $ so  that $ x=0 $ is an inflow boundary and $ x=1 $ is  an outflow boundary. We  equally distribute  $ M+1 $ points $\{x_j = (j+1/2)h: j=0,1,\dots,M\}$ in the computational domain $(0,1)$, as shown in Figure \ref{fig:kreiss}. We use $u_j$ to denote the value of $u$ at $x=x_j$ and suppress the index for the time levels. Obviously
at the inflow boundary, the solution value at the ghost point $ x_{-1} $ is required in order  to perform a second-order finite difference at $ x_0 $. For this purpose, a polynomial is constructed in the region around the inflow boundary by using point-wise values $ u_{-1} $, $ u_0 $ and $ u_1 $,
\begin{equation}
L(x) = g_{-1}(x)u_{-1} + g_0(x)u_0 + g_1(x)u_1,
\end{equation}
from which we want to find  the value $ u_{-1} $. The Lagrangian interpolation tells that
\begin{equation}
\begin{array}{l}
g_{-1}(x) = \dfr{(x-x_0)(x-x_1)}{2h^2}, \\[3mm]  g_0(x) = \dfr{(x-x_{-1})(x-x_1)}{-h^2}, \\[3mm]  g_1(x) = \dfr{(x-x_{-1})(x-x_0)}{2h^2}.
\end{array}
\end{equation}
Then $ u_{-1} $ can be obtained by solving the linear equation $  u(0,t) = L(0) $ where  $ u(0,t) = g(t)$. At the outflow boundary $ x=1 $, we simply use the extrapolation
\begin{equation}
  u_{M+1} = 2u_M - u_{M-1}
\end{equation}
 to obtain the value $ u_{M+1} $ since the signal goes out of the computational domain at this end.
 \vspace{0.2cm}

 The extension to high order is highly nontrivial. Let's review the approach developed in \cite{Du-Li-B-2018}. The same as other multi-stage methods (e.g. in \cite{Don-1995}), the current two-stage fourth order method  needs careful treatment at the intermediate stage for the approximation in order to preserve the accuracy.   There are two key points: The construction at ghosts points and the approximation of boundary conditions.  The inverse LW approach is applied here \cite{Tan-Shu-ILW-2010}, but the current treatment is much simpler and easier to be implemented since no higher derivatives need calculating.

 \vspace{0.2cm}

 \subsection{Ghost values}

 We still consider IBVP \eqref{eq:IBVP}. Then the values  $ \bar{u}_{-1} $, $ \bar{u}_{-2} $, $ \Delta u_{-1} $ and $\Delta u_{-2}$, {defined over  $I_{-1}=(-h,0)$ and $I_{-2}=(-2h,-h)$}, ghost cells,   are needed in the reconstruction procedure for the values indexed by $ j = 0 $ and $ j = 1 $.  To obtain the values mentioned above, a cubic polynomial,
\begin{equation}\label{eq:poly-5}
 p(x) = \al_3 x^3 +\al_2 x^2 +\al_1 x + \al_0,
\end{equation}
is constructed over $ I_{-2} \cup I_{-1} \cup I_0 \cup I_1 = (-2h,2h) $ to interpolate the solution $ u(x,t^n) $
such that
\begin{equation}\label{eq:condition-5}
 \dfr{1}{h}\displaystyle\int_{I_{i}}p(x)dx = \bar{u}_{i}, \ \ \ \ i=-2, -1, 0, 1.
\end{equation}
With the constraints \eqref{eq:condition-5} into \eqref{eq:poly-5}, we determine the coefficients  $ \al_0 $, $ \al_1 $, $ \al_2 $ and $ \al_3 $  as
\begin{equation}\label{eq:poly-5-coeff-0}
\begin{array}{ll}
\alpha_3 = \dfr{\bar{u}_1 - 3\bar{u}_0 + 3\bar{u}_{-1} - \bar{u}_{-2}}{6h^3}, \ \ &\alpha_2 = \dfr{\bar{u}_1 - \bar{u}_0 - \bar{u}_{-1} + \bar{u}_{-2}}{4h^2},\\[3mm]
\alpha_1 = \dfr{-\bar{u}_1 + 15\bar{u}_0 - 15\bar{u}_{-1} + \bar{u}_{-2}}{12h}, \ \ &\alpha_0 = \dfr{-\bar{u}_1 + 7\bar{u}_0  + 7\bar{u}_{-1}  -\bar{u}_{-2}}{12},
\end{array}
\end{equation}
in which $ \bar{u}_{-1} $ and $ \bar{u}_{-2} $ are yet to be determined and they are obtained by evaluating $ p(0) $ and $ p^\prime(0) $ at the boundary $x=0$, 
\begin{equation}\label{eq:boundary-interpolation}
\begin{array}{l}
  p(0) = \dfr{1}{12}(-\bar{u}_{1} + 7\bar{u}_{0} + 7\bar{u}_{-1} - \bar{u}_{-2}) = g(t)+\mathcal{O}(h^4),\\[2mm]
  p^\prime(0) = \dfr{1}{12h}(-\bar{u}_{1}+15\bar{u}_{0}-15\bar{u}_{-1}+\bar{u}_{-2})= -f^\prime(g(t))^{-1} \ g^{\prime}(t)+\mathcal{O}(h^3).
\end{array}
\end{equation}
Solving \eqref{eq:boundary-interpolation} in terms of $\bar{u}_{-1}$ and $\bar{u}_{-2}$ yields (by ignoring high order terms)
\begin{equation}\label{eq:ILW-result-u}
\begin{array}{l}
  \bar{u}_{-1} = \dfr 14 (-6g + 6~h~ f'(g)^{-1} \ g^{\prime} + 11\bar{u}_{0} - \bar{u}_{1}),\\[3mm]
  \bar{u}_{-2} = \dfr 14 (-90g + 42~h~ f'(g)^{-1} \ g^{\prime} + 105\bar{u}_{0} - 11\bar{u}_{1}).
\end{array}
\end{equation}
Substituting \eqref{eq:ILW-result-u} into \eqref{eq:poly-5-coeff-0}, in turn,  gives us the explicit expressions of $ \alpha_i$, $i=0,\dots, 3$, and then the expression of $ p(x) $. Therefore we have (by ignoring high order terms)
\begin{equation}\label{eq:ILW-result-ux}
\bga{l}
  \Delta u_{-1} = \dfr{p(0) - p(-h)}h\\[2mm]
\ \ \ \ \ \ \ \ = \dfr 1{8h} (66g - 34~h~ f'(g)^{-1} \ g^{\prime} - 73\bar{u}_{0} + 7\bar{u}_{1}),\\[2mm]
  \Delta u_{-2} = \dfr{p(-h) - p(-2h)}h\\[2mm]
\ \ \ \ \ \ \ \ = \dfr 1{8h} (294g - 118~h~ f'(g)^{-1} \ g^{\prime} - 331\bar{u}_{0} + 37\bar{u}_{1}).
\eda
\end{equation}
 Thus \eqref{eq:ILW-result-u} and \eqref{eq:ILW-result-ux} together provide the  values in the ghost cells $I_{-1}$ and $I_{-2}$. Not that in \eqref{eq:boundary-interpolation} the inverse LW approach is used,
\vspace{2mm}

As there are discontinuities close to the inflow boundary, a WENO-type stencil selecting procedure can be applied. Assume that there is a discontinuity in either $ I_0 $ or $ I_1 $, we  shorten the stencil cell by cell. Denote the stencils by
\begin{equation}
\bga{l}
S^{(2)} = \{I_{-2}, I_{-1}, I_0, I_1\},\ \
S^{(1)} = \{I_{-2}, I_{-1}, I_0\},\ \
S^{(0)} = \{I_{-2}, I_{-1}\}.
\eda
\label{stencil}
\end{equation}
Denote by $ p^{(r)} (x)$  the interpolation polynomial  on  $ S^{(r)} $, $r=0,1,2$, just as the polynomial $ p(x) $ constructed before.   Then define
\begin{equation}
\bga{ll}
  \bar{ u}_{-1}^{(r)} = \dfr{1}{h}\displaystyle\int_{I_{-1}}p^{(r)}(x)dx, & \bar{ u}_{-2}^{(r)} = \dfr{1}{h}\displaystyle\int_{I_{-2}}p^{(r)}(x)dx,\\[3mm]
  \Delta u_{-1}^{(r)} = \dfr 1h (p^{(r)}(0) - p^{(r)}(-h)), & \Delta u_{-2}^{(r)} = \dfr 1h (p^{(r)}(-h) - p^{(r)}(-2h)).
\eda
\end{equation}
The expressions of $ \bar{ u}_{-1}^{(r)} $, $ \bar{ u}_{-2}^{(r)} $, $ \Delta  u^{(r)}_{-1} $ and $ \Delta  u^{(r)}_{-2} $ for $ r=0,1,2 $ will be listed in \ref{app:low-order}.

The smoothness indicators are defined in the same way as for the classical WENO interpolation, and the values are given
\begin{equation}\label{eq:weno-ilw-un1}
\bga{ll}
  \bar{ u}_{-1} = \d\sum_{r=0}^2\omega^{(r)}\bar{ u}_{-1}^{(r)},& \bar{ u}_{-2} = \d\sum_{r=0}^2\omega^{(r)}\bar{ u}_{-2}^{(r)},\\
  \Delta  u_{-1} = \d\sum_{r=0}^2\omega^{(r)}\Delta  u_{-1}^{(r)},& \Delta  u_{-2} = \d\sum_{r=0}^2\omega^{(r)}\Delta  u_{-2}^{(r)},
\eda
\end{equation}
where the linear weights of each stencil are
\begin{equation}
\bga{l}
\label{eq:weight}
\al^{(r)} = \dfr{d^{(r)}}{{(\varepsilon+\beta^{(r)})}^2},\ \ \ \
\omega^{(r)} = \dfr{\al^{(r)}}{\sum_{l=0}^2{\al^{(l)}}}\\[3mm]
 d^{(0)} = h^2,\ \ d^{(1)} = h,\ \ d^{(2)} = 1 - d^{(0)} - d^{(1)}.
\eda
\end{equation}

\subsection{Inflow boundary condition treatment at intermediate stages}\label{sec:inter-stage}
The same as in other multi-stage temporal discretization  \cite{BD-ACC,BD-ACC-2}, the direct use of exact boundary conditions at intermediate stages in the process of multi-stage approaches will cause the lose of the numerical accuracy.  In order to offset such a defect,
our strategy is made as follows.  We first focus on  the leftmost control volume $I_0$ and write out the solution advancing formula,
\begin{equation}
\begin{array}{rl}
\bar{u}^{n+1}_0 & =\bar u_0^n -\dfr{ k }{h} \left[f_{\frac 12}^{4th} -f_{-\frac 12}^{4th}\right]\\[3mm]
 &= \bar u_0^n -\dfr{1}{h} \Big\{ k\left[f(u_{\frac 12}^n) -f(u_{-\frac 12}^n)\right]\\[3mm]
 & \ \ \ \ \ \  \ \   +\dfr{ k^2 }{6}\left[f^\prime(u^{n}_{\frac 12})(\dfr{\pt u}{\pt t})^{n}_{\frac 12} -  f^\prime( u^{n}_{-\frac 12})(\dfr{\pt u}{\pt t})^{n}_{-\frac 12}\right]\\[3mm]
&\ \ \ \ \ \ \ \ +\left.\dfr{ k^2 }{3}\left[f^\prime( u^{n+\frac 12}_{\frac 12})(\dfr{\pt u}{\pt t})^{n+\frac 12}_{\frac 12} - f^\prime(u^{n+\frac 12}_{-\frac 12})(\dfr{\pt u}{\pt t})^{n+\frac 12}_{-\frac 12}\right]  \right\}.
 \end{array}
 \end{equation}
Using the governing equation \eqref{eq:IBVP} to replace the temporal derivatives by the corresponding spatial ones, we obtain
\begin{equation}\label{eq:1d-al2-dire}
\begin{array}{rl}
\bar{u}^{n+1}_0&= \bar u_0^n -\dfr{1}{h} \Big\{ k\left[f(u_{\frac 12}^n) - f(u_{-\frac 12}^n)\right]\\[3mm]
 & \ \ \ \ \ \  \ \ -\dfr{ k^2 }{6}\left[(f^\prime(u^{n}_{\frac 12}))^2(\dfr{\pt u}{\pt x})^{n}_{\frac 12} - (f^\prime(u^{n}_{-\frac 12}))^2(\dfr{\pt u}{\pt x})^{n}_{-\frac 12}\right] \\[3mm]
& \ \ \ \ \  \ \ -\dfr{ k^2 }{3}\left[(f^\prime(u^{n+\frac 12}_{\frac 12}))^2(\dfr{\pt u}{\pt x})^{n+\frac 12}_{\frac 12} - (f^\prime(u^{n+\frac 12}_{-\frac 12}))^2(\dfr{\pt u}{\pt x})^{n+\frac 12}_{-\frac 12}\right] \Big\}.
\end{array}
\end{equation}
The difficulty results from the presence of $ ({\pt u/\pt x})^{n+\frac 12}_{\frac 12} $ and $ ({\pt u/\pt x})^{n+\frac 12}_{-\frac 12} $  evaluated at the intermediate stage $t=t^{n+\frac 12}$.
In order to restore the fourth-order accuracy of the two-stage fourth-order scheme,   we use
\begin{equation}\label{eq:derv-recon-dire-m}
\begin{array}{l}
  {\left(\dfr{\pt u}{\pt x}\right)^{n+\frac 12}_{-\frac 12}} = -(f^\prime(g(t^{n+\frac{1}{2}})))^{-1}(g^{\prime})^{n+\frac{1}{2}},\\[3mm]
  {\left(\dfr{\pt u}{\pt x}\right)^{n+\frac 12}_{\frac 12}} = \dfr{1}{48h}\left[-49\bar{u}^{n+\frac{1}{2}}_{0} + 59\bar{u}^{n+\frac{1}{2}}_{1} - 4\bar{u}^{n+\frac{1}{2}}_{2} \right. \\[3mm]
   \ \ \ \ \ \ \ \ \ \ \ \ \ \ \ \ \ \ \ \ \ \ \ \left. - 6g^{n+\frac{1}{2}} + 6h(f^\prime(g(t^{n+\frac{1}{2}})))^{-1}(g^{\prime})^{n+\frac{1}{2}}\right],
\end{array}
\end{equation}
where the exact boundary values $g(t^{n+\frac 12})$ and $g^\prime(t^{n+\frac 12})$   are replaced by
 \begin{equation}\label{eq:inter-boundary-value}
\begin{array}{l}
  g^{n+\frac{1}{2}} = g(t^{n+\frac{1}{2}}) -\dfr{ k ^3}{48}g^{\prime\prime\prime}(t^{n+\frac 12}),\\[2mm]
  {(g^{\prime})}^{n+\frac{1}{2}} = g^\prime(t^{n+\frac{1}{2}}).\\
\end{array}
\end{equation}
The detailed analysis is given in \cite{Du-Li-B-2018}.

\subsection{Outflow boundary condition}
We set   $ x_{M+\frac 12}=1 $ as an outflow boundary, at which  no boundary condition is prescribed theoretically. Numerically,  we have to set  required values $ \bar{u}_{M+1} $, $ \bar{u}_{M+2} $, $ \Delta u_{M+1} $ and $ \Delta u_{M+2} $ in ghost cells.
Since the signal propagates out of the computational domain through the boundary $ x=1 $, the extrapolation can be  used to construct the data in the  ghost cells  $I_{M+1}$ and $I_{M+2}$.
 A cubic polynomial is constructed  in order to achieve the fourth-order accuracy,
 \begin{equation}\label{eq:extra-4-poly}
\bga{l}
  q(x) = \dfr{\bar{u}_{M}-3\bar{u}_{M-1}+3\bar{u}_{M-2}-\bar{ u}_{M-3}}{6h^3}(x-1)^3\\[3mm]
\ \ \ \ \ \ \ \ + \dfr{5\bar{u}_{M}-13\bar{u}_{M-1}+11\bar{u}_{M-2}-3\bar{u}_{M-3}}{4h^2}(x-1)^2\\[3mm]
\ \ \ \ \ \ \ \  + \dfr{35\bar{u}_{M}-69\bar{u}_{M-1}+45\bar{u}_{M-2}-11\bar{u}_{M-3}}{12h}(x-1)\\[3mm]
\ \ \ \ \ \ \ \  + \dfr{25\bar{u}_{M}-23\bar{u}_{M-1}+13\bar{u}_{M-2}-3\bar{u}_{M-3}}{12}.
\eda
\end{equation}
This gives the values
\begin{equation}\label{eq:extra-4-res}
\bga{l}
  \bar{u}_{M+1} = 4\bar{u}_{M} - 6\bar{u}_{M-1} + 4\bar{u}_{M-2} - \bar{u}_{M-3},\\[3mm]
  \bar{u}_{M+2} = 10\bar{u}_{M} - 20\bar{u}_{M-1} + 15\bar{u}_{M-2} - 4\bar{u}_{M-3},\\[3mm]
  \Delta u_{M+1} = \dfr{26\bar{u}_{M} - 57\bar{u}_{M-1} + 42\bar{u}_{M-2} - 11\bar{u}_{M-3}}{6},\\[3mm]
  \Delta u_{M+2} = \dfr{47\bar{u}_{M} - 114\bar{u}_{M-1} + 93\bar{u}_{M-2} - 26\bar{u}_{M-3}}{6}.
\eda
\end{equation}
If there is a discontinuity in either $ I_{M-3} $, $ I_{M-2} $, $ I_{M-1} $ or $ I_{M} $, a WENO-type stencil selection can be applied.
\vspace{0.2cm}

\subsection{Hyperbolic systems}

At moment, the boundary treatment for  systems of hyperbolic balance laws is basically achieved through the diagonalization process. Then we distinguish various cases such as  the solid boundary condition, inflow and outflow boundary conditions for practical applications.   Details can be found in \cite{Du-Li-B-2018}.

\section{Computational Performance}

In our series of papers, we have demonstrated the performance of current temporal-spatially coupled algorithms through many challenging benchmark problems, particularly in \cite{Pan-Li-2017} using the GKS solver.  Here I would like to give some remarks in terms of computational efficiency, robustness and fidelity.

\subsection{Computational Efficiency}

Computational efficiency is always an important issue for practical engineering problems.  We have tested and compared the efficiency with the popular WENO algorithm in \cite{Xu-Li-2017}   and with DG in  \cite{Cheng-Li-2018}.

Specifically, in \cite{Xu-Li-2017} we evaluate the computational costs of the WENO-type reconstruction and the flux evaluation quantitatively.
The time for each reconstruction is denoted by $T_R$, the
time for second-order gas-kinetic solver is $T_{2nd}$, and the time
for third-order flux solver is $T_{3rd}$. According to the data provided
in  \cite[Table 1, Page 203]{Xu-Li-2017}, we can estimate the time used for the computations of flux and
reconstruction with the following relations,
\begin{align*}
T_R+2T_{2nd}=0.84287s,\\
T_R+2T_{3rd}=1.38178s,\\
2T_R+12T_{2nd}=2.20566s,
\end{align*}
where the estimation is based on the characteristic variable reconstruction and each flux is shared by two cells.
Thus, the time for reconstruction is $T_R=0.71289s$, the time for
second-order gas-kinetic flux solver is $T_{2nd}=0.06499s$ and  the time
for third-order gas-kinetic flux solver is $T_{3rd}=0.33445s$. For
classical fourth-order Runge-Kutta schemes, the computational time for four
spatial reconstruction alone will become much higher than the fourth-order
gas-kinetic scheme  for the update of each cell averaged values
\begin{align*}
4T_R = 2.85156s > 2T_R+ 12 T_{2nd} = 2.20566s.
\end{align*}
Similar estimation can be done for the conservative variables reconstruction.
Even without counting on the cost of the flux evaluation in the traditional fourth-order Runge-Kutta method,
such as those commonly used with the Lax-Friedrichs flux,
the current fourth-order time stepping method  is still more efficient than the classical methods.

\begin{figure}[!htb]
  \begin{center}
    \includegraphics[width=7cm]{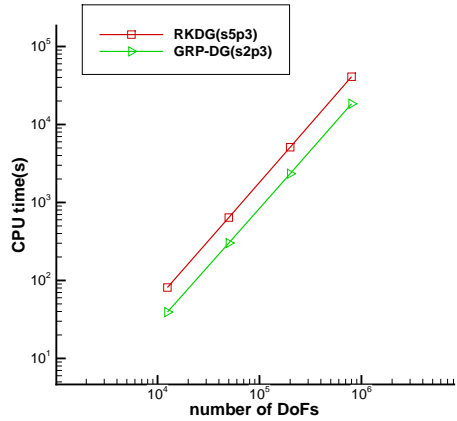}
  \end{center}
  \caption{Comparison of CPU time(s) between RKDG($s_5p_3$) method and GRP-DG($s_2p_3$) method for the shock vortex interaction problem.}\label{fig:ex5-52}
\end{figure}

The efficiency  is mainly attributed to the half of reconstruction steps compared to that for the same order of other line methods.  This is further verified in the framework of DG methods \cite{Cheng-Li-2018}. In Table \ref{tab:ex5-5} and Figure \ref{fig:ex5-52} through simulating shock-vortex interaction problem,    demonstrating  that nearly $55\%$ CPU time can be saved using
 the  GRP-DG(s2p3) method  compared to the same order SSP RKDG(s5p3)
method. This result meets the expectation well as compared to the RKDG(s5p3) method which needs
five stages of evaluating DoFs and performing reconstruction to achieve fourth order, while the GRP-DG(s2p3) method only takes two stages to provide totally comparable results.

\begin{table}[!htb]
\centering \caption{Comparison of CPU time(s) between RKDG and GRP-DG methods for shock vortex interaction problem.}\label{tab:ex5-5}
\begin{tabular}{c|c|c|c|c}
  Methods             & $50\times25$ & $100\times50$ & $200\times100$ & $400\times200$ \\
\hline
  GRP-DG($s_2p_3$)    &   39.3    &  301.6  &  2339.4  & 18443.4 \\
\hline
  RKDG($s_5p_3$)      &   81.2    &  640.7  &  5109.3  & 40999.6  \\
\end{tabular}
\end{table}

\vspace{0.2cm}

\subsection{Robustness} The robustness is always an important indicator for a practical numerical method.  In the framework of multistage multi-derivative algorithm, the strong stability preserving (SSP) property was taken over to show the stability \cite{Seal-2018}. However,  SSP seems not work when the Lax-Wendorff type flow solvers are taken as the building block.  Therefore new stability framework is worth exploring in the future.
\vspace{0.2cm}

Nevertheless, in practice, the current two-stage fourth order accurate algorithm has the same stability as the second order version: the Courant number  is taken above $0.5$ except  extreme cases such as the large density ration problem. Empirically, the current ``$2\odot 2$" algorithm is more robust than other multi-stage methods.

\subsection{Fidelity} 
\begin{figure}[!htp]
\centering
\includegraphics[width=11cm]{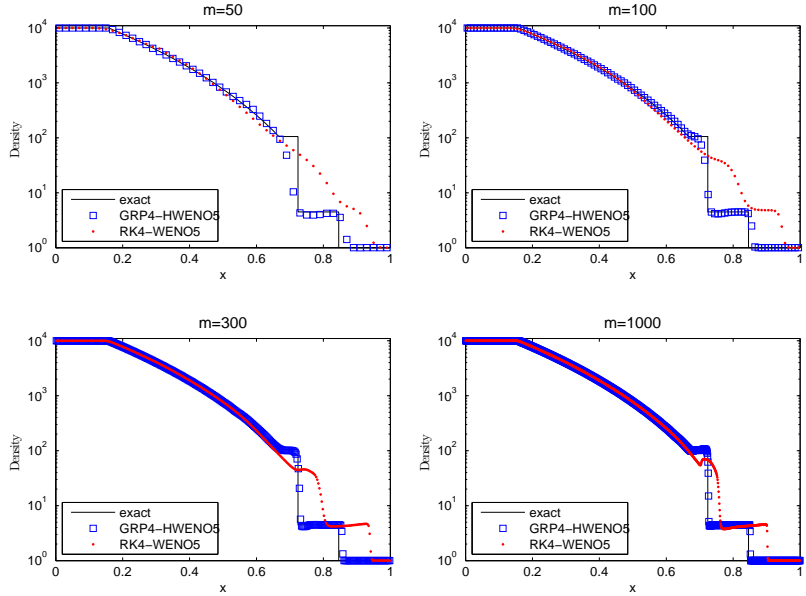}
  \caption[small]{The comparison of the density profile for the large pressure ratio problem in Example 6. The schemes are GRP4-HWENO5 (squares) and RK4-WENO5 (dots) with $ m $ cells.The solid lines are the exact solution.}
\label{fig:ratio}
\end{figure}

In the community of CFD,  the {\em fidelity} is termed for a numerical simulation of very complex problems  using a specific algorithm. Since there is no  reliable  mathematical theory in general supporting the current CFD simulation, the verification of high fidelity  appears very valuable.  We pursue such studies in the whole process of the current algorithm. For example, we resolve the associated GRP analytically and use the GRP solution for Hermite-type data reconstruction. In \cite{Li-Wang-2017} we elaborate the so-called large density ratio problem \cite{Liu-Tang-2006} using the GRP solver. When  the current ``$2\odot 2$" algorithm is adopted, quite few grids are needed to obtain satisfactory results, as shown in Figure \ref{fig:ratio}.

Also readers are recommended to test the benchmark problems in \cite{Pan-Li-2017}, for which all simulations are made in the ``$2\odot 2"$ framework using the GKS solver.

\section{Prospective Discussions }

It is natural to require the temporal-spatial coupling of a numerical method when simulating compressible fluid flows,  for which the GRP solver and the GKS solver are reviewed briefly as the representatives of Lax-Wendroff type flow solvers.
The direct embedding into any numerical frameworks,  such as the finite volume/DG framework, already results in  favorable second order numerical schemes. Interested readers can refer to papers by Jiequan Li and his collaborators for GRP methods (www.ams.org/mathscinet, scholar.google.com or researchgate.net).
\vspace{0.2cm}

As to the ``$2\odot 2$" algorithm itself, it is just at the beginning stage, and many issues are awaiting for our study. Below are some immediate doable problems.
\begin{enumerate}
\item[$P_1$] What is a good framework for stability analysis?

\item[$P_2$]  It is valuable to compare and develop multi-stage two-derivative algorithms with arbitrary order of accuracy.

\item[$P_3$] Develop implicit ``$2\odot 2$" algorithm with various applications such as detonation simulation.

\item[$P_4$]  Apply this algorithm for the simulation of turbulence flows and other engineering problems.

\end{enumerate}
You are welcome to join this new branch of high order numerical methods for CFD.

\appendix
\section{The interpolation results in Subsection \ref{sec:inflow}}
This appendix is dedicated to list the interpolation results  in Section \ref{sec:inflow}. Recall that we assume $ x=0 $ and $ x=1 $ are the inflow and outflow boundaries for the IBVP \eqref{eq:IBVP} of the one-dimensional scalar conservation laws, respectively. The stencils are denoted in \eqref{stencil}.

\

\subsection{Cell averages and cell differences}\label{app:low-order}
The reconstructed average of $ u $ in $ I_{-1} $ and $ I_{-2} $ on those stencils are:
\begin{equation}
\bga{l}
      \bar{u}_{-1}^{(2)} = \dfr 14 (-6g + 6~h~ f^\prime(g)^{-1} \ g^{\prime} + 11\bar{u}_{0} - \bar{u}_{1}),\\[2mm]
      \bar{u}_{-1}^{(1)} = h~ f^\prime(g)^{-1} \ g^{\prime} + \bar{u}_{0},\\[2mm]
      \bar{u}_{-1}^{(0)} = g + \dfr 12~ h~ f^\prime(g)^{-1} \ g^{\prime},\\[4mm]
      \bar{u}_{-2}^{(2)} = \dfr 14 (-90g + 42~h~ f^\prime(g)^{-1} \ g^{\prime} + 105\bar{u}_{0} - 11\bar{u}_{1}),\\[2mm]
      \bar{u}_{-2}^{(1)} = -6gh + 5~f^\prime(g)^{-1} \ g^{\prime} + 7\bar{u}_{0},\\[2mm]
      \bar{u}_{-2}^{(0)} = g + \dfr 32~ h~ f^\prime(g)^{-1} \ g^{\prime}.
\eda
\end{equation}
The reconstructed $ x $-difference of $ u $ in $ I_{-1} $ and $ I_{-2} $ on those stencils are:
\begin{equation}
\bga{l}
      \Delta u_{-1}^{(2)} = \dfr 1{8h} (66g - 34~h~ f^\prime(g)^{-1} \ g^{\prime} - 73\bar{u}_{0} + 7\bar{u}_{1}),\\[2mm]
      \Delta u_{-1}^{(1)} = \dfr 1{2h} (6g - 5~h~ f^\prime(g)^{-1} \ g^{\prime} - 6\bar{u}_{0}),\\[2mm]
      \Delta u_{-1}^{(0)} = -f^\prime(g)^{-1} \ g^{\prime},\\[4mm]

      \Delta u_{-2}^{(2)} = \dfr 1{8h} (294g - 118~h~ f^\prime(g)^{-1} \ g^{\prime} - 331\bar{u}_{0} + 37\bar{u}_{1}),\\[2mm]
      \Delta u_{-2}^{(1)} = \dfr 1{2h} (18g - 11~h~ f^\prime(g)^{-1} \ g^{\prime} - 18\bar{u}_{0}),\\[2mm]
      \Delta u_{-2}^{(0)} = -f^\prime(g)^{-1} \ g^{\prime}.
\eda
\end{equation}

\subsection{Smoothness indicators}\label{app:SI}
The smoothness indicators on these stencils are listed as follows,
\begin{equation}
\bga{l}
\beta^{(2)} = \dfr 1{80}\Big[    66516 g^2 +     9444 (h f^\prime(g)^{-1}g^\prime)^2 - 56348  f^\prime(g)^{-1}g^\prime h \bar{ u}_0\\[2mm]
\ \ \ \ \ \ \ \ \ \ \ \ \ \ \ \ \ + 85929 \bar{ u}_0^2 + 6644  f^\prime(g)^{-1}g^\prime h \bar{ u}_1 -     20694 \bar{ u}_0 \bar{ u}_1 +     1281 \bar{ u}_1^2\\[2mm]
\ \ \ \ \ \ \ \ \ \ \ \ \ \ \ \ \ + 12 g (    4142  f^\prime(g)^{-1}g^\prime h -    12597 \bar{ u}_0 +     1511 \bar{ u}_1)\Big],\\[3mm]

\beta^{(1)} = 48 g^2 + 54 g h  f^\prime(g)^{-1}g^\prime + 16 (h f^\prime(g)^{-1}g^\prime)^2\\[2mm]
\ \ \ \ \ \ \ \ \ \ \ \ \ - 96 g \bar{ u}_0 + 48 \bar{ u}_0^2 - 54 h f^\prime(g)^{-1}g^\prime \bar{ u}_0,\\[3mm]

\beta^{(0)} = (h f^\prime(g)^{-1}g^\prime)^2.
\eda
\end{equation}

\begin{acknowledgements}
This work is supported by NSFC (nos. 11771054, 91852207) and Foundation of LCP.  This  review paper is based on the joint works with Matania Ben-Artzi, Jian Cheng, Zhifang Du, Xin Lei, Liang Pan,  Jin Qi, Yue Wang and K. Xu, to whom the author expresses his deep thanks.
\end{acknowledgements}



\end{document}